\documentclass[11pt]{article}

\usepackage[usenames,dvipsnames]{color}

\usepackage{amssymb, mathrsfs}
\usepackage{amsthm}
\usepackage{amsmath,color}
\usepackage{graphicx, tikz-cd}
\usepackage{yfonts}
\usepackage{amsbsy}
\usepackage{hyperref}
\usepackage{multicol}
\usepackage[margin=1in]{geometry}
\usepackage{listings}
\usepackage{comment} 
\usepackage[utf8]{inputenc}
\usepackage{float}
\usepackage{breqn}
\usepackage{cleveref}
\makeatletter
\newcommand{\address}[1]{\gdef\@address{#1}}
\newcommand{\email}[1]{\gdef\@email{\url{#1}}}
\newcommand{\@endstuff}{\par\vspace{\baselineskip}\noindent\small
\begin{tabular}{@{}l}\scshape\@address\\\textit{E-mail address:} \@email\end{tabular}}
\AtEndDocument{\@endstuff}
\makeatother
\address{Department of Mathematics, Rice University}
\email{yandi.wu@rice.edu}

\theoremstyle{plain}
\newtheorem{theorem}{Theorem}[section]
\newtheorem{corollary}{Corollary}[theorem]
\newtheorem{lemma}[theorem]{Lemma}

\newtheorem{proposition}[theorem]{Proposition}

\theoremstyle{definition}
\newtheorem{assumption}[theorem]{Assumption}
\newtheorem{definition}[theorem]{Definition}
\newtheorem{remark}[theorem]{Remark}

\usepackage[backend=biber,style=ieee-alphabetic,sorting=nyvt]{biblatex}
\allowdisplaybreaks
\begin{document}

\title{Marked Length Spectrum Rigidity for Surface Amalgams} 
\date{}
\author{Yandi Wu} 

\maketitle 

\begin{abstract}
    In this paper, we show that simple, thick negatively curved two-dimensional P-manifolds, a large class of surface amalgams, are marked length spectrum rigid. That is, if two piecewise negatively curved Riemannian metrics (satisfying certain smoothness conditions) on a simple, thick two-dimensional P-manifold assign the same lengths to all closed geodesics, then they differ by an isometry up to isotopy. Our main theorem is a natural generalization of Croke and Otal's celebrated results about marked length spectrum rigidity of negatively curved surfaces.  
\end{abstract}

\section{Introduction}

Mostow's famous rigidity theorem states that for closed hyperbolic manifolds of dimension greater than two, the metric is completely determined by the fundamental group up to isotopy. Mostow's rigidity theorem does not hold for hyperbolic surfaces or negatively curved Riemannian metrics without constant sectional curvature. For these cases, in order to conclude two metrics are equivalent up to isotopy, further restrictions are needed, such as requiring the surfaces to have the same \textit{marked length spectra}, defined as follows:  
 
\begin{definition}[Marked length spectrum] The \textit{marked length spectrum of a metric space $(M, g)$} is the class function 
$$MLS: \pi_1(M) \rightarrow \mathbb{R}^{+}, [a] \mapsto \inf\limits_{{\gamma \in [\alpha]}} \ell_g(\gamma)$$ 
which assigns to each free homotopy class $[\alpha] \in \pi_1(M)$ the infimum of lengths in the class.
\end{definition}

In particular, if $g$ is negatively curved or locally CAT(-1), each homotopy class of curves has a unique geodesic representative, so the marked length spectrum assigns a length to each closed geodesic in $(M, g)$. We say that $(M, g_0)$ and $(M, g_1)$ have the same marked length spectrum if for every $[
\alpha] \in \pi_1(M)$, $MLS_0([\alpha]) = MLS_1([\alpha])$. A class of metrics $\mathcal{M}$ is \textit{marked length spectrum rigid} if whenever $(M, g_0)$ and $(M, g_1) \in \mathcal{M}$ have the same marked length spectrum, there exists an isometry $\varphi: (M, g_0) \rightarrow (M, g_1)$ that is isotopic to the identity. \\

In this paper, we will study a class of objects, defined below, that are natural generalizations of surfaces. The following definition is adapted from Definition 2.3 of \cite{LF07}. 

\begin{definition}[Negatively curved two-dimensional P-manifold]\label{pmnfld} A compact metric space $X$ is a \textit{negatively curved two-dimensional P-manifold} if there exists a closed subset $Y \subset X$ (the \textit{gluing curves} of $X$) that satisfies the following: 

\begin{enumerate}
    \item Each connected component of $Y$ is homeomorphic to $S^1$;
    \item The closure of each connected component of $X - Y$ is homeomorphic to a compact surface with boundary endowed with a negatively curved (Riemannian) metric, and the homeomorphism takes the component of $X - Y$ to the interior of a surface with boundary. We will call each $\overline{X - Y}$ a \textit{chamber} in $X$; 
    \item There exists a negatively-curved metric on each chamber which coincides with the original metric. 
\end{enumerate}
\end{definition}

If $Y$ forms a totally geodesic subspace of $X$ consisting of disjoint simple closed curves, we say that $X$ is \textit{simple}. If each connected component of $Y$ (gluing curve) is attached to at least three distinct boundary components of chambers, then we say $X$ is \textit{thick} (note this definition differs slightly from the one given in \cite{LF07} but follows the one given in \cite{LF2} and \cite{LF3}). Like Lafont in \cite{LF07}, we will only be considering simple, thick two-dimensional P-manifolds. Doing so ensures $X$ is locally CAT(-1) (in other words, its universal cover, $\widetilde{X}$, is CAT$(-1)$) and guarantees some useful properties, such as one pointed out in Lemma \ref{disks}. Lafont has defined higher dimensional P-manifolds as well (see \cite{LF2}); we focus on the two-dimensional ones, which consist of surfaces with boundary glued together along their boundary components. Throughout this paper, all P-manifolds will be assumed to be two-dimensional.

\begin{figure}
    \centering
    \includegraphics[width=\textwidth]{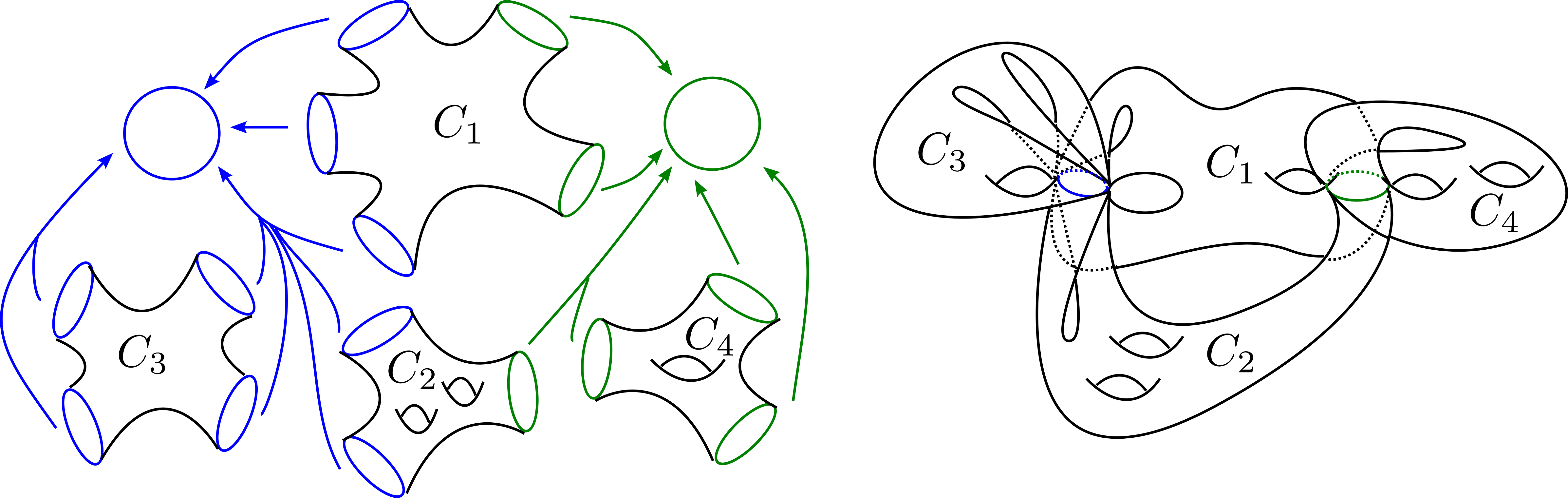}
    \caption{An example of a simple, thick (two-dimensional) P-manifold with four chambers.}
    \label{fig:pmnfld}
\end{figure}

We will equip a simple, thick negatively-curved P-manifold $X$ with a metric in a class we denote as $\mathcal{M}_{\leq}$, following the notation from \cite{con}. For the precise definition of $\mathcal{M}_{\leq}$, we refer the reader to section \ref{Mleq}, but roughly speaking, metrics in $\mathcal{M}_{\leq}$ are locally CAT(-1) piecewise Riemannian metrics with an additional condition that limits pathological behavior around the gluing curves of $X$. We now state the main result of the paper: 

\begin{theorem}\label{main} Suppose $(X, g_1)$ and $(X, g_2)$ are simple, thick negatively-curved P-manifolds equipped with metrics $g_1, g_2 \in \mathcal{M}_{\leq}$. Furthermore, suppose $(X, g_1)$ and $(X, g_2)$ have the same marked length spectrum. Then there is an isometry $\phi: (X, g_1) \rightarrow (X, g_2)$ that is isotopic to the identity. 

\end{theorem}

The marked length spectrum rigidity problem has a long history. Fricke and Klein showed the class of closed hyperbolic surfaces is marked length spectrum rigid (see \cite{FK65}). Burns and Katok then conjectured that closed negatively curved manifolds of all dimensions are marked length spectrum rigid in \cite{bk}. A major breakthrough occurred when Croke and Otal independently proved the conjecture in dimension two (see \cite{Otal90}, \cite{croke}). Croke and Otal's results led to a large quantity of generalizations of the marked length spectrum rigidity problem for surfaces, some (but far from all) of which are listed here. Croke, Fathi, and Feldman extended Croke's result to the case of non-positively curved metrics in \cite{cff}. In another direction, Hersonksy and Paulin extended Otal's result to the case of negatively curved metrics with finitely many cone singularities in \cite{hp}. Duchin, Leininger, and Rafi proved marked length spectrum rigidity for closed surfaces endowed with flat metrics with finitely many cone singularities arising from quadratic differentials in \cite{dlr}, a result extended by Bankovic and Leininger in \cite{bl}, who proved the result for all flat metrics with finitely many cone singularities. Finally, by combining previous methods, Constantine proved marked length spectrum rigidity for nonpositively curved metrics on surfaces with finitely many cone singularities in \cite{con2}. We remark that marked length spectrum rigidity is unlikely to hold for surfaces with boundary, although Guillarmou and Mazzucchelli have shown a weaker form of rigidity, \textit{marked boundary rigidity}, for a large family of negatively curved metrics on surfaces with strictly convex boundary (see \cite{gm}).

While the marked length spectrum rigidity problem has been well-studied in the case of surfaces, for higher dimensions, the conjecture remains largely open. Hamenst\"adt showed marked length spectrum rigidity for rank one locally symmetric closed manifolds of dimension greater than 2 (see \cite{ham}). There is also a local marked length spectrum rigidity result for closed, compact manifolds of all dimensions with Anosov geodesic flow due to Guillarmou and Lefeuvre (see \cite{gl}). For metric spaces that are not manifolds, the question also remains largely unstudied, although there are a few notable results. Work by Culler and Morgan (see \cite{cm}) and Alperin and Bass (see \cite{ab}) leads to a marked length spectrum rigidity result for Culler-Vogtmann Outer Space, built to study the group Out($\mathbb{F}_n$) in analogy to the relationship between the Teichm\"uller space of $S$, $\mathscr{T}(S)$, and the mapping class group, $\text{Mod}(S)$. Generalizing the work from \cite{cm} and \cite{ab}, Constantine and Lafont prove a version of marked length spectrum rigidity of compact, non-contractible one-dimensional spaces in \cite{cl2}. In \cite{con}, they show that certain compact quotients of Fuchsian buildings, including those with piecewise hyperbolic metrics, are marked length spectrum rigid. \\

\noindent\textbf{Outline of the Paper.} In Section \ref{prelim}, we define precisely the class of metrics $\mathcal{M}_{\leq}$ and review some fundamental facts about CAT(-1) spaces. Furthermore, we detail some separability properties of hyperbolic groups developed by Haglund and Wise which we will use in \Cref{generalcase}. In Section \ref{easy}, we patch together isometries constructed in \cite{Otal90} to prove the theorem for P-manifolds with the property that each chamber can be included into a closed surface. Since there is no well-defined unit tangent bundle, proving ergodicity of the geodesic flow map, a key component of Otal's original proof, requires some heavy machinery included in the Appendix. Finally, in Section \ref{generalcase}, we prove the general case by constructing finite-sheeted covers to reduce to the base case examined in Section \ref{easy}. The argument relies on separability properties of the fundamental groups of simple, negatively-curved surface amalgams, which can be realized as non-positively curved cube complexes. In particular, we show that such groups are \textit{QCERF}, which allows one to promote certain immersions to embeddings in finite-sheeted covers. \\

\noindent \textbf{Acknowledgements.} First and foremost, I am deeply indebted to Caglar Uyanik for many helpful and illuminating discussions concerning the project, and for his edits and comments on previous drafts of this paper. Furthermore, I would like to thank Chris Leininger for his sharp and invaluable insight, and for contributing ideas towards the general case of the argument. I would also like to thank Tullia Dymarz for many useful conversations and insights, and help with proof checking. Many thanks also to Jean-Francois Lafont for generously suggesting the problem as a sequel to his work on P-manifolds, and for very helpful discussions during a visit to the Ohio State University. Thanks also goes to David Constantine for discussing the project with me in depth. Thanks goes to Alex Hof for suggesting the idea behind the proof of Lemma \ref{crossratio}. Additionally, I sincerely thank the anonymous referee for carefully reading and proofreading the paper and pointing out a gap in my original argument. Finally, I gratefully acknowledge support from an NSF RTG grant DMS-2230900.
\section{Preliminaries} 
\label{prelim}
\subsection{The class of metrics $\mathcal{M}_{\leq}$}
\label{Mleq}
We first precisely define the class of metrics $\mathcal{M}_{\leq}$ under consideration in Theorem \ref{main} and point out some important properties. We say $g \in \mathcal{M}_{\leq}$ if $g$ satisfies the following properties: 

\begin{enumerate}
    \item Each chamber of $C \subset X$ is equipped with a negatively-curved Riemannian metric with sectional curvature bounded above by $-1$ so that $C$ has geodesic boundary components; 
    \item The restrictions of $g$ to the chambers of $X$ are ``compatible" in the sense that if two boundary components $b_1$ and $b_2$ of two (possibly the same) chambers $C_1$ and $C_2$ are both attached to a gluing curve $\gamma \subset X$, then the gluing maps $b_1 \hookrightarrow \gamma$ and $b_2 \hookrightarrow \gamma$ are isometries (in particular, we do not allow circle maps of degree two or greater like those explored in \cite{hst});
    \item For any two boundary components $b_1 \in C_1$ and $b_2 \in C_2$, the restriction of $g$ to $N_{b_1} \bigcup\limits_{b_1 \sim b_2} N_{b_2}$ is a negatively curved smooth Riemannian metric with sectional curvature bounded above by $-1$, where $N_{b_1}$ and $N_{b_2}$ are $\epsilon$-neighborhoods around $b_1$ and $b_2$ respectively for some $\epsilon > 0$. 
    
\end{enumerate}

We impose the third condition to ensure that we can exploit previous marked length spectrum rigidity results for surfaces which in particular require Riemannian negatively curved metrics with at most a finite number of cone singularities. We remark that Theorem \ref{main} still holds if one introduces a finite number of cone singularities (points with cone angles greater than $2\pi$) to each chamber; instead of using Theorem 1 from \cite{Otal90} in the proof of Lemma \ref{chamberisom}, one would instead use Theorem C from \cite{hp} or Corollary 2 from \cite{con2}. Furthermore, we suspect Condition 3 could be eliminated if one were to carefully modify Otal's proof to allow for metrics that are Riemannian outside a singular set of gluing geodesics, but more work would need to be done. 

Finally, we point out that a surface amalgam is \textit{non-positively curved (NPC)} if each chamber admits a non-positively curved Riemannian metric. Note that the only difference between negatively curved and NPC surface amalgams is the latter may have chambers which are cylinders. For the sake of completeness, we will briefly allow NPC surface amalgams in some auxilliary lemmas in \Cref{generalcase}.  

We now discuss some properties of $\mathcal{M}_{\leq}$ that will be useful in the proof of Theorem \ref{main}. 

\subsubsection{Properties of $\mathcal{M}_{\leq}$}

First, we remark that metrics in $\mathcal{M}_{\leq}$ are locally CAT(-1):

\begin{remark}\label{cat-1}
If $(X, g)$ is a negatively-curved P-manifold where $g \in \mathcal{M}_{\leq}$, then $(X, g)$ is locally CAT(-1).  
\end{remark}

Indeed, suppose $X$ is equipped with a metric $g \in \mathcal{M}_{\leq}$ and $C \subset X$ is a chamber in $X$. Recall a generalization of the Cartan-Hadamard Theorem which states that a smooth Riemannian manifold $M$ has sectional curvature $\leq \kappa$ if and only if $M$ is locally CAT($\kappa$) (see \cite{bh} Theorem 1A.6). As a result, the restriction of $g$ to $C$ is locally CAT(-1) since $C$ is a endowed with a negatively curved metric with sectional curvature bounded above by $-1$.  If $\kappa \in \mathbb{R}$ and $X_1$ and $X_2$ are locally CAT($\kappa$) spaces glued isometrically along a convex, complete metric subspace $A \subset X_1 \cap X_2$, then $X_1 \sqcup_A X_2$ is locally CAT($\kappa$) (see Theorem 2.11.1 in \cite{bh}). As a result, since we have endowed each chamber in a P-manifold $X$ with a locally CAT(-1) metric, $(X, g)$ will also be locally CAT(-1), which shows the remark.

The locally CAT(-1) property of metrics in $\mathcal{M}_{\leq}$ will prove useful in the proof of the base case of the main result in the Section \ref{easy}. Note that CAT(-1) spaces are also Gromov ($\delta$-)hyperbolic, so $(X, g)$ where $g \in \mathcal{M}_{\leq}$ is also Gromov hyperbolic, another useful property we will exploit in Section \ref{easy}. 

\begin{remark} \label{gromhyp}
Another way to show that $(X, g)$ where $g \in \mathcal{M}_{\leq}$ is Gromov hyperbolic is to observe that any two-dimensional P-manifold (a metric space satisfying only Conditions 1 and 2 from Definition \ref{pmnfld}) has a fundamental group that is an amalgamated product of surface groups and free groups over cyclic subgroups (without any $\mathbb{Z}^2$ subgroups), which is Gromov hyperbolic by the Bestvina-Feighn Combination Theorem (see \cite{bf}). 
\end{remark}

\subsection{The Universal Cover of a Simple, Thick P-manifold}

Next, we describe the universal cover of a simple, thick P-manifold, which we will be working with extensively in Section \ref{easy}. Roughly speaking, $\widetilde{X}$ will look like an infinite collection of totally geodesic subspaces of disks (homeomorphic to $\mathbb{H}^2$) glued together in a tree-like fashion. Following Lafont, we will call disks in $\widetilde{X}$ \textit{apartments}. Throughout this paper, also following Lafont, we will also call polygonal lifts of chambers in the universal cover a \textit{chamber}. Following Lafont, we will call geodesics in $\widetilde{X}$ that project to gluing curves under the covering map \textit{branching geodesics}. 

\begin{figure}[h!]
 \begin{centering}
\includegraphics[width=0.8\textwidth]{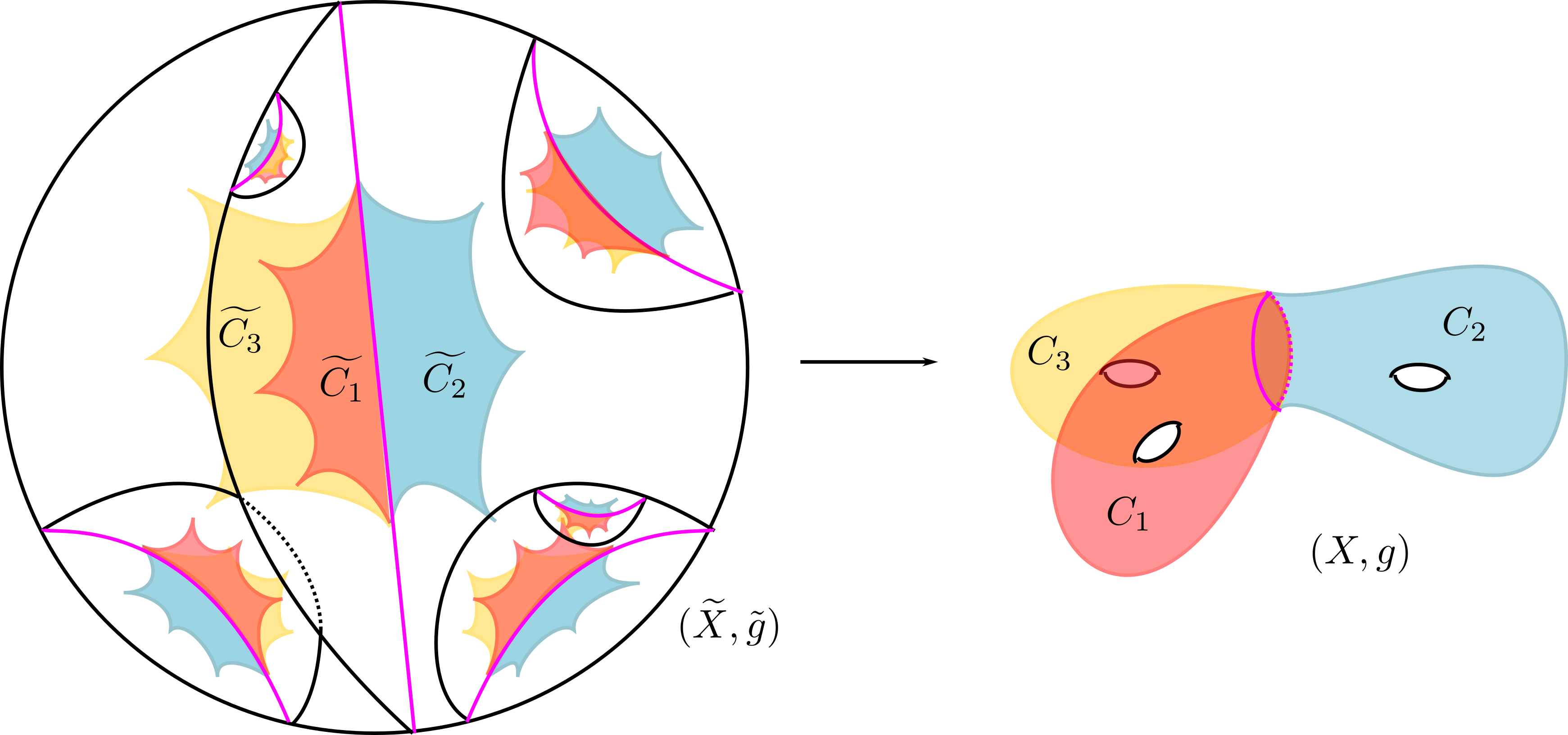}
    \caption{The universal cover of a simple, thick P-manifold $(X, g)$.}
    \label{fig:cover}
\end{centering} 
\end{figure}

While in general, $\widetilde{X}$ is very hard to visualize, the case where each chamber can be included into a closed surface is much easier to describe. We will focus on describing $\widetilde{X}$ in this special case since our proof strategy relies on reducing to it. Each gluing geodesic $\gamma$ will lift to a countably infinite collection of branching geodesics in $\widetilde{X}$, with a countably infinite collection in each apartment containing lifts of chambers that intersect $\gamma$. Each of these lifts of $\gamma$ are adjacent to $n$ half planes, where $n$ is the number of boundary components attached to $\gamma$. Each closed surface will lift to an infinite collection of apartments tiled by lifts of chambers that are embedded in the closed surface. See Figure \ref{fig:cover} for an illustration of $\widetilde{X}$. 

\begin{remark}
Some readers may notice similarities between $\widetilde{X}$ and a two-dimensional hyperbolic building. Even the terminology of P-manifolds borrows heavily from that of buildings; both have ``chambers," ``apartments" and sets of branching geodesics, which are called \textit{walls} in building terminology. This allows us to borrow some definitions, such as intersection numbers, from \cite{con}. However, for buildings, boundaries of chambers are not totally geodesic in the space, so their walls are not totally geodesic. The proof methods from \cite{con} are different from those seen in this paper in part due to this fact. Furthermore, all the chambers in buildings are assumed to be isometric, which is not the case in $\widetilde{X}$. See Section 1.4.3 of \cite{lfthesis} for a more complete list of differences between buildings and P-manifolds. 
\end{remark}

\subsection{Visual Boundaries}
\label{boundary}

In this section, we will assume that $(\widetilde{X}, \widetilde{g})$ is a proper (closed balls are compact) metric space that is either CAT(-1) or Gromov hyperbolic; the definitions are the same for both. For more details and background, we refer the reader to \cite{bk2}. We say that two geodesic rays in $\widetilde{X}$, $\alpha_1: \mathbb{R}_{\geq 0} \rightarrow \widetilde{X}$, and $\alpha_2: \mathbb{R}_{\geq 0} \rightarrow \widetilde{X}$ are \textit{asymptotic} if they lie within bounded distance of one another; in other words, there exists some finite $M \geq 0$ such that for all $t \in \mathbb{R}_{\geq 0}$, $d(\alpha_1(t), \alpha_2(t)) \leq M$. To define $\partial^{\infty}(\widetilde{X}, \widetilde{g})$, fix any basepoint $x_0 \in (\widetilde{X}, \widetilde{g})$ and consider the set of all geodesic rays originating from $x_0$. Two geodesic rays based at $x_0$ are equivalent if they are asymptotic. 

\begin{definition}[Visual Boundary]\label{boundary} We define equivalence classes of such geodesic rays based at $x_0$ as the \textit{visual boundary of $\widetilde{X}$}, $\partial^{\infty}(\widetilde{X}, \widetilde{g})$.
\end{definition} 

When $\widetilde{g}$ is CAT(-1), there is a unique geodesic in $(\widetilde{X}, \widetilde{g})$ between any two points in $\partial^{\infty}(\widetilde{X}, \widetilde{g})$ (see Proposition 1.4 of \cite{bh}), so the space of oriented geodesics in $\widetilde{X}$ is identified with $\partial^{\infty}(\widetilde{X}, \widetilde{g}) \times \partial^{\infty}(\widetilde{X}, \widetilde{g}) \setminus \Delta$, where $\Delta$ indicates the diagonal.   \\

We define the \text{cross ratio} as follows: 

\begin{definition}[Cross ratio and M\"{o}bius maps] \label{crm} If $(\widetilde{X}, \widetilde{g})$ has boundary $\partial^{\infty}(\widetilde{X}, \widetilde{g})$, one can define the \textit{cross ratio} of a four-tuple of distinct boundary points $(\xi, \xi', \eta, \eta') \in (\partial^{\infty}(\widetilde{X}))^4$: 
\begin{equation}
\label{cr}
[\xi\xi'\eta\eta'] = \lim\limits_{(a,a',b,b') \rightarrow (\xi,\xi',\eta,\eta')} \frac{1}{2}\big(\widetilde{g}(a, b) + \widetilde{g}(a', b') - \widetilde{g}(a, b') - \widetilde{g}(a', b)\big).
\end{equation} 
See Figure \ref{fig:crossratio} for an example in $\mathbb{H}^2$). 

 \begin{figure}[h!]
 \begin{centering}
\includegraphics[width=0.3\textwidth]{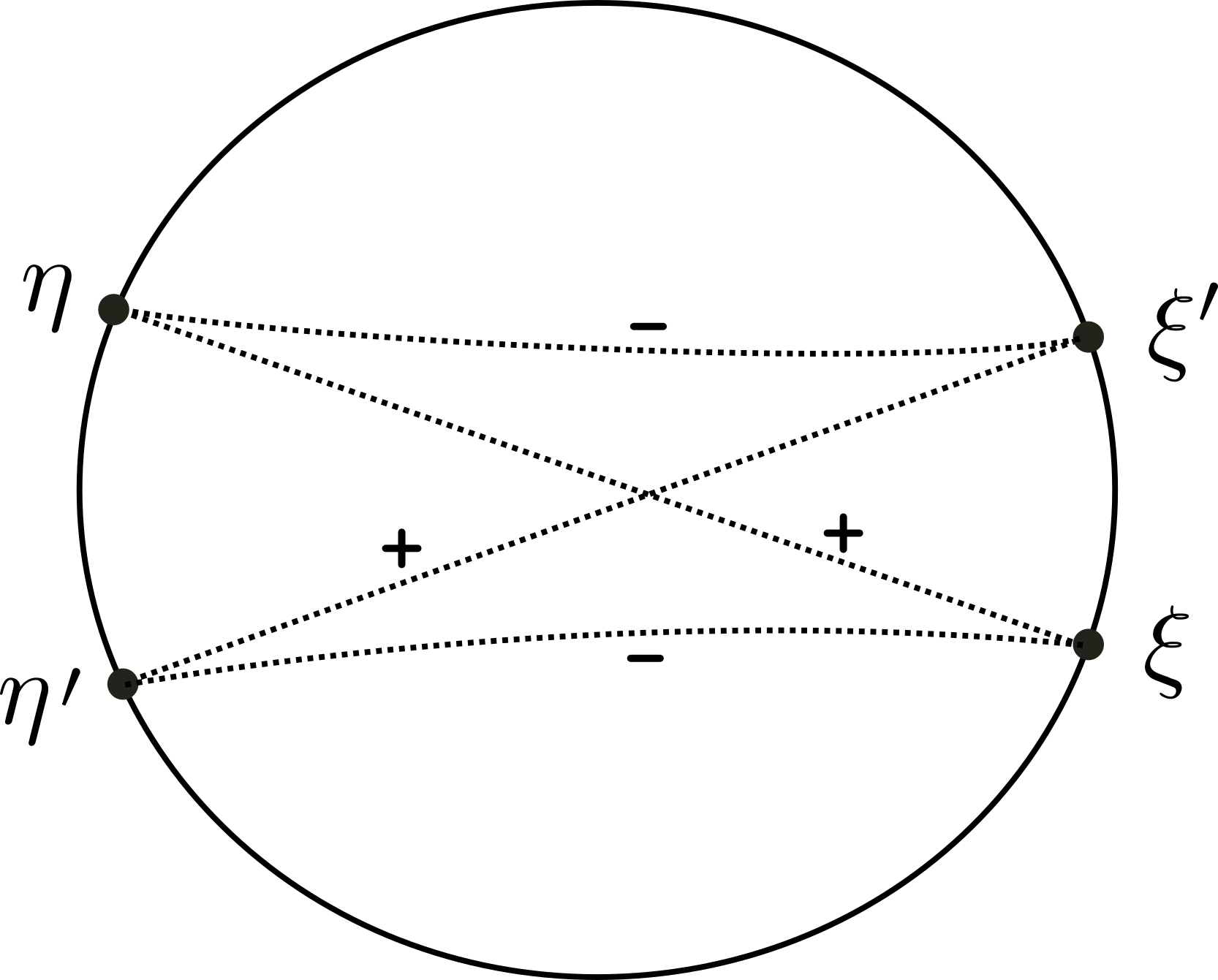}
    \caption{An illustration of the calculation of $[\xi\xi'\eta\eta']$ in $\mathbb{H}^2$.}
    \label{fig:crossratio}
\end{centering} 
\end{figure}

A map $\partial^{\infty}f: \partial^{\infty}(\widetilde{X}, \widetilde{g_1}) \rightarrow \partial^{\infty}(\widetilde{X}, \widetilde{g_2})$ is \textit{M\"{o}bius} if it preserves the cross ratio (i.e. $[\xi\xi'\eta\eta'] = [f(\xi)f(\xi')f(\eta)f(\eta')]$). 
\end{definition}

From the definition of the cross-ratio, it follows that: 

\begin{equation}\label{crsymm} 
[abcd] = [cdab] \text{ and } -[abdc] = [abcd] = -[bacd].
\end{equation}

We are now ready to define the Gromov Product: 

\begin{definition}[Gromov Product]
\label{gromovproduct} Let $x \in (\widetilde{X}, \widetilde{g})$ and $a, b \in \partial^{\infty}(\widetilde{X})$. Then the \textit{Gromov Product} of $a$ and $b$ with respect to a basepoint $x \in \widetilde{X}$ is:
$$ 
\langle a, b \rangle_{x} = \lim\limits_{(a_i, b_i) \rightarrow (a, b)} \frac{1}{2}\big(\widetilde{g}(a_i, x) + \widetilde{g}(x, b_i) - \widetilde{g}(a_i, b_i)\big).
$$

\end{definition} 

The Gromov product measures how long two geodesics travel close together in $(\widetilde{X}, \widetilde{g})$. Indeed, if $a_i, b_i, x \in \widetilde{X}$ are three arbitrary points in a $\delta$-hyperbolic metric space, then $$\langle a_i, b_i \rangle_x = \frac{1}{2}\big(\widetilde{g}(a_i, x) + \widetilde{g}(x, b_i) - \widetilde{g}(a_i, b_i)\big)$$ approximates within $2\delta$ the distance between $x$ and the geodesic segment $[a_i, b_i]$. As a result, the Gromov product provides a convenient way to topologize $\partial^{\infty}(\widetilde{X}, \widetilde{g})$. Indeed, we can define a countable neighborhood base for any $p \in \partial^{\infty}(\widetilde{X})$ and $d \in \mathbb{N}$ as follows: 

\begin{equation}
    \label{nbhd}
    B(p, d) = \{ q \in \partial^{\infty}(\widetilde{X}) \vert \langle p, q\rangle_x  > d\}.
\end{equation}

Using the Gromov Product, we can also endow $\partial^{\infty}(\widetilde{X}, \widetilde{g})$ with a \textit{visual metric} that induces the aforementioned topology on $\partial^{\infty}(\widetilde{X}, \widetilde{g})$: 

\begin{definition}[Visual metric] \label{visualmetric}
Let $(\widetilde{X}, \widetilde{g})$ be a proper CAT(-1) space. Suppose $a, b \in \partial^{\infty}(\widetilde{X}, \widetilde{g})$. Then for some fixed basepoint $x \in (\widetilde{X}, \widetilde{g})$, we can assign $\partial^{\infty}(\widetilde{X}, \widetilde{g})$ a \textit{visual metric}:

\begin{equation}\widetilde{g}_{\infty, x}(a, b) = \begin{cases} e^{-\langle a, b \rangle_x} \text{ if $a \neq b$}\\
0 \text{ otherwise}.
\end{cases}
\end{equation}
\end{definition}

Note that if $x' \in (\widetilde{X}, \widetilde{g})$ is a different basepoint, then we have that for $A = e^{2\delta}e^{\widetilde{g}(x, x')} > 1$:

$$A^{-1}\widetilde{g}_{\infty, x}(a, b) \leq \widetilde{g}_{\infty, x'}(a, b) \leq A\widetilde{g}_{\infty, x}(a, b).$$

\subsubsection{The Visual Boundary of a P-manifold}
\label{visual}

We now present some useful properties of the visual boundary of a simple, thick P-manifold $(\widetilde{X}, \widetilde{g})$. As before, we will assume all the chambers in $X$ are negatively curved. 

Bass-Serre theory provides a useful characterization of points on $\partial^{\infty}(\widetilde{X})$. Consider a bipartite \textit{graph of groups decomposition} $\mathcal{G}$ of $G = \pi_1(X)$, where there is a vertex group $\pi_1(C)$ for each chamber $C \subset X$ and $\langle \gamma \rangle \cong \mathbb{Z}$ for each gluing curve $\gamma \subset X$. Furthermore, there is a cyclic edge group between $\langle \gamma \rangle$ and $\pi_1(C)$ for each boundary component of $C$ that is attached to $\gamma$. Thus, $\mathcal{G} = V_0 \sqcup V_1$ can be partitioned into collection of vertices labeled with $\langle \gamma \rangle$, which we will call $V_0$, and those labeled with $\pi_1(C)$, which we will call $V_1$ (see \cite{dst} for more details on and examples of this construction). The \textit{Bass-Serre tree} of $G$ is a tree $T$ on which $G$ acts minimally (i.e. there is no proper invariant subtree of $T$) and without inversions on edges with quotient $\mathcal{G} = T / G$. 

To better describe the Bass-Serre tree $T$, we briefly summarize the main points from Section 4.1 of \cite{malone}, which uses Bowditch's construction to characterize boundary points for \textit{geometric amalgams of free groups} (e.g. $\pi_1(X)$). Given a point in $x \in \partial^{\infty}(\widetilde{X})$, Bowditch characterizes $x$ by $\text{val}(x)$, the number of topological ends of $\partial^{\infty}(\widetilde{X}) \setminus \{x\}$. In particular, for a simple, thick P-manifold with negatively-curved chambers, $\text{val}(x) = n \geq 3$ if $x$ is a branching geodesic attached to $n$ distinct boundary components and $\text{val}(x) = 2$ otherwise. Bowditch defines separate equivalence classes for $M(2) = \{x \in \partial^{\infty}(\widetilde{X}) : \text{val}(x) = 2\}$ and $M(3+) = \{ x \in \partial^{\infty}(\widetilde{X}) : \text{val}(x) \geq 3\}$:

\begin{enumerate}
\item \textit{For points in $M(3+)$}: We say $x \approx y$ if either $x = y$ or the number of connected components of $\partial^{\infty}(\widetilde{X}) \setminus \{x, y\}$ is equal to $\text{val}(x)$ and $\text{val}(y)$. 
\item \textit{For points in $M(2)$}: We say $x \sim y$ if $x = y$ or the number of connected components of $\partial^{\infty}(\widetilde{X}) \setminus \{x, y\}$ is equal to $2$.  
\end{enumerate} 

In both cases, equivalent pairs of points form a \text{cut pair} of $\partial^{\infty}(\widetilde{X})$. The Bass-Serre tree is a bipartite infinite-valence tree consisting of vertices labeled by equivalence classes in $M(3+)$ and $M(2)$, which represent conjugacy classes of vertices in $V_0$ and $V_1$ respectively.  Two vertices $\Lambda_e$ and $\Lambda_v \in T$ corresponding to a $\approx$-class in $M(3+)$ and $\sim$-class in $M(2)$ respectively are connected by an edge if any two distinct points $x, y \in \Lambda_v$ lie on the same connected component of $\partial^{\infty}(\widetilde{X}) \setminus \{a, b\}$, where $a, b$ are the two distinct points in $\Lambda_e$ corresponding to endpoints of some branching geodesic. Bowditch thus formulates a trichotomy of points in $\partial^{\infty}(\widetilde{X})$: 

\begin{proposition}\label{bassserre}
[Proposition 1.3 of \cite{bowditch} and Lemma 10 of \cite{kk}] If $\pi_1(X)$ acts on Bass-Serre tree $T$, then each point $x \in \partial^{\infty}(\widetilde{X})$ corresponds to exactly one of the following equivalence classes of points in $\partial^{\infty}(\widetilde{X})$: 

\begin{enumerate}
    \item A point in $\Lambda_e$ corresponding to some conjugate of $V_0 \in V(\mathcal{G})$;  
    \item A point in $\Lambda_v$ corresponding to some conjugate of $V_1 \in V(\mathcal{G})$; 
    \item A point of $\partial^{\infty}(T)$, with a unique $x$ for each point in $\partial^{\infty}(T)$.
\end{enumerate}
\end{proposition}

In other words, $\approx$-classes in $M(3+)$ correspond to vertices in category (1) in Proposition \ref{bassserre} while vertices in category (2) correspond to $\sim$-classes of $M(2)$. 

Using Proposition \ref{bassserre}, we can topologically characterize points in $\partial^{\infty}(\widetilde{X})$. As mentioned earlier, points in category (1) correspond to endpoints of branching geodesics, lifts of gluing curves in $\partial^{\infty}(\widetilde{X})$. Equivalence classes of points corresponding to a vertex in category (2) consist of points $x$ and $y$ such that the unique geodesic $(x, y)$ between $x$ and $y$ does not intersect any branching geodesics. Indeed, if $(x, y)$ were to intersect a branching geodesic, $\partial^{\infty}(\widetilde{X}) \setminus \{x, y\}$ would still be connected. Topologically, one can check that this means $x$ and $y$ in fact lie on the boundary of the same connected component of $p^{-1}(C) \setminus \{\mathcal{BG}\}$, where $C$ is some chamber in $X$, $p: \widetilde{X} \rightarrow X$ is the universal covering map, and $\mathcal{BG}$ is the full set of branching geodesics in $\widetilde{X}$. A vertex $\Lambda_e$ is connected to a vertex $\Lambda_v$ if the branching geodesic corresponding to $\Lambda_e$ is a boundary component of the connected component of $p^{-1}(C) \setminus \{\mathcal{BG}\}$ associated with $\Lambda_v$. From this information, we see that traveling along an edge in $T$ is equivalent to crossing a branching geodesic in $\widetilde{X}$. Thus, the unique points in category (3) correspond to geodesic rays that intersect infinitely many branching geodesics in $\widetilde{X}$. 

We now point out that points on $\partial^{\infty}(\widetilde{X})$ lie on boundaries of embedded disks, a fact stated in the proof of Lemma 3.1 in \cite{LF07}. For the convenience of the reader, we provide a proof as well. 

\begin{lemma}
\label{disks}
Given a geodesic $\gamma \in \widetilde{X}$ with endpoints $p, q \in \partial^{\infty}(\widetilde{X})$, there exists an isometrically embedded apartment containing $\gamma$. That is, $p$ and $q$ lie in an embedded $S^1 \subset \partial^{\infty}(\widetilde{X})$. As a consequence, the boundary $\partial^{\infty}(\widetilde{X})$ is path-connected.
\end{lemma}

\begin{proof} First, we consider the case where $\gamma$ is a gluing geodesic. If $p(\gamma)$, the gluing geodesic lifting to $\gamma$ in $\widetilde{X}$, is a closed geodesic in a closed surface, then there is automatically an apartment containing $\gamma$. Otherwise, consider two distinct boundary components $b_1 \subset C_1$ and $b_2 \subset C_2$ attached to $p(\gamma)$; note that $b_1$ and $b_2$ exist due to the thickness hypothesis and $C_1$ and $C_2$ are not required to be distinct chambers. In $\widetilde{X}$, there is a set of lifts of $C_1$, which we will call $K_1$, that forms a totally geodesic subset of a half-disk containing $\gamma$. Similarly, there is a set of lifts of $C_2$ forming a totally geodesic subset $K_2$ (disjoint from $K_1$) of a half-disk containing $\gamma$. Note that even though $K_1 \cup K_2$ is only a subset of a disk containing $\gamma$, we can ``fill out"  $K_1 \cup K_2$ to obtain an apartment $A$ embedded in $\widetilde{X}$. Indeed; if there is part of a disk missing from $K_1 \cup K_2$, its boundary must necessarily be a branching geodesic $\gamma'$ in $\widetilde{X}$. By the thickness hypothesis, there is some collection of polygons $\mathcal{P}$ disjoint from $K_1$ and adjacent to $\gamma'$ that project to a chamber $C$ adjacent to a gluing geodesic $p(\gamma')$ in $X$. We can ``continue" $A$ by attaching a subset of a half-plane disjoint from $K_1$ and tiled by copies the polygons in $\mathcal{P}$. We then iterate, attaching collections of lifts of chambers and eventually limiting to a half-plane $H_1 \supset K_1$ embedded in $\partial^{\infty}(\widetilde{X})$. Apply the same procedure to obtain a half-disk $H_2$ containing $K_2$. Note that this construction may be counterintuitive because $A = H_1 \cup H_2$ does not project under the covering map to any surfaces in $X$ (see Figure \ref{fig:fillin}). 

The case where $\gamma$ is not a gluing geodesic is similar. Note that in this case, $\gamma$ will pass through a sequence of branching geodesics $\{\gamma_n\}_{n \in \mathbb{N}}$. For each $\gamma_n$ that $\gamma$ passes through, $\gamma$ will also pass through a collection $\mathcal{P}_n$ of lifts of $C_n$, a chamber attached to $p(\gamma_n)$. Note $\mathcal{P}_n$ can be chosen to be a totally geodesic subspace of a disk, $K_n$, containing both $\gamma_n$ and $\gamma_{n + 1}$; similar to the procedure from before, we can ``fill in" $K_n$ to a connected section of a disk bordered by $\gamma_n$ and $\gamma_{n + 1}$, which we will call $H_n$. Iterate the procedure to obtain an apartment $A = \bigcup\limits_{n \in \mathbb{N}} H_n$. 
For a proof of why this implies $\partial^{\infty}(\widetilde{X})$ is path-connected, we refer the reader to \cite{LF07}.  

\begin{figure}[h!]
    \centering
    \includegraphics[width=0.8\textwidth]{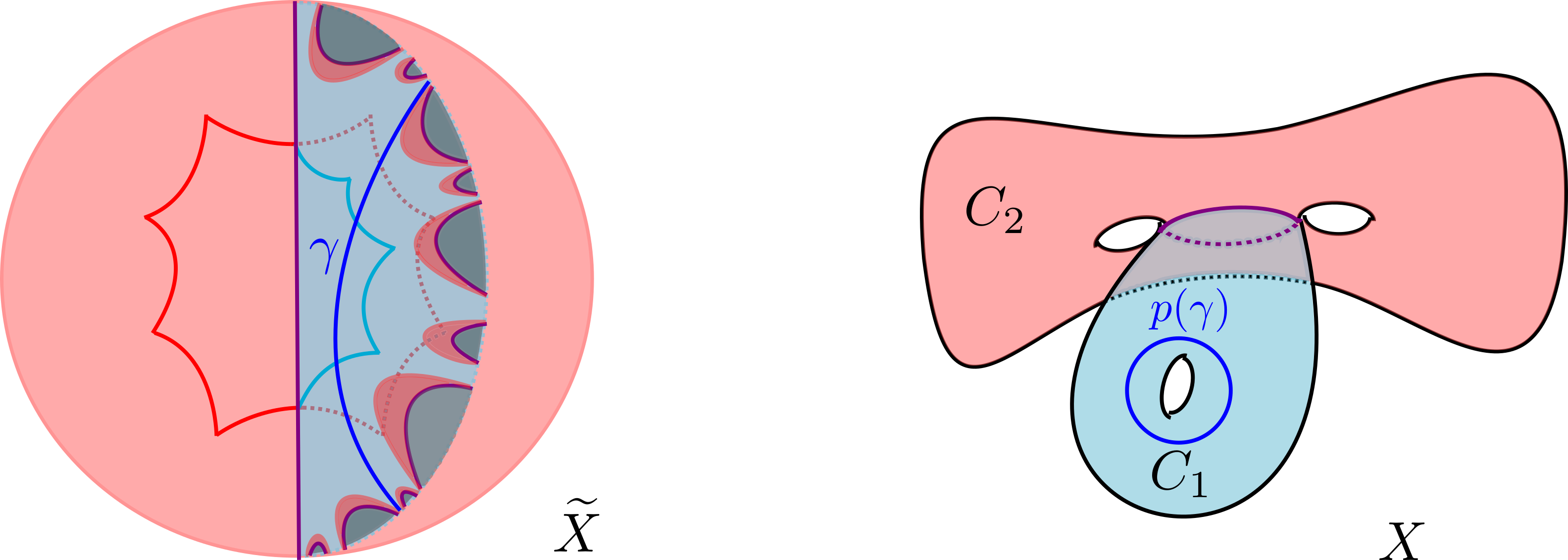}
    \caption{An example of how to choose an apartment containing $\gamma$. In this example, $p(\gamma)$ is completely contained in the torus with one boundary component $C_1$, which lifts to a totally geodesic subspace of a half-disk. We can fill in the missing portions of the half-disk (shaded in gray) with portions of disks that are copies of the universal cover of $C_2$. While this is an especially straightforward case of the filling in procedure, the construction works in general.}
    \label{fig:fillin}
\end{figure}
    
\end{proof}

Lemma \ref{disks} allows us to define intervals in $\partial^{\infty}(\widetilde{X})$: 

\begin{definition}[Interval on $\partial^{\infty}(\widetilde{X})$]
    \label{interval} 
Suppose $X$ is a simple, thick negatively-curved P-manifold. An \textit{interval} on $\partial^{\infty}(\widetilde{X})$ with endpoints $a, b \in \partial^{\infty}(\widetilde{X})$ is an interval on the boundary of an apartment containing the unique geodesic in $\widetilde{X}$ with endpoints $a$ and $b$. 
\end{definition}

\subsection{Subgroup Separability, QCERF, and NPC Cube Complexes}
\label{sec:qcerf}

We now briefly introduce key tools for the arguments in \Cref{generalcase} and give a preview of some claims we will prove. 

\begin{definition}[Separability of a Subgroup] Let $G$ be a group and let $H \leq G$ be a subgroup. We say that $H$ is \textit{separable} in $G$ if for any $g \in G \setminus H$, there exists a finite-index subgroup $K < G$ which contains $H$ but $g \notin K$.    
\end{definition}

If the trivial subgroup $H = \{1\}$ is separable in $G$, we say that $G$ is \textit{residually finite}. Separability properties have attracted considerable interest in recent decades due to an associated powerful topological property which allows one to promote certain immersions to embeddings: 

\begin{theorem}[Scott, \cite{scott}, Lemma 1.3]\label{thm:scott} Let $G$ be a finitely generated group and $X$ a CW complex such that $\pi_1(X) = G$. Let $H \leq G$ and $Y \rightarrow X$ a cover corresponding to $H$. Then $H$ is separable if and only if for every finite subcomplex $K \subseteq Y$, there exists an intermediate finite-sheeted cover $Y \rightarrow Z \rightarrow X$ such that $K$ embeds as a subcomplex of $Z$. 
\end{theorem}

A group $G$ is \textit{subgroup separable}, or, for historical reasons, LERF, if every finitely-generated subgroup of $G$ is separable. In general, it is unknown which fundamental groups of simple, thick, negatively curved surface amalgams are subgroup separable, although results due to Gitik in \cite{gitik} provide partial results. For example, the fundamental groups of the surface amalgams from \Cref{fig:cover} and \Cref{fig:fillin} are subgroup separable. 

When subgroup separability of a group $G$ is unknown, one can sometimes show that well-behaved subgroups of $G$ are separable. A subspace $Y$ of a geodesic metric space $X$ is \textit{quasiconvex} if there exists a 
$K \geq 0$ such that, for all $y_1, y_2 \in Y$ and all $x \in [y_1, y_2]$, $d(x, Y) \leq K$. Note that if $K = 0$, $Y$ is convex. A group $H$ acting on a geodesic metric space $X$ is \textit{quasiconvex} if the orbit $Hx$ is a quasiconvex subspace of $X$ for some (any) $x \in X$. We point out a stronger notion: if $H$ acts cocompactly on and stabilizes $C$, a convex subset of $X$, then $H$ is \textit{convex}. A group $G$ is \textit{quasiconvex subgroup separable (QCERF)} if every quasiconvex subgroup of $G$ is separable. 

A particularly well-behaved class of groups are those which act properly and cocompactly on \textit{non-positively curved (NPC) cube complexes}, also known as CAT(0) cube complexes. Recall that a \textit{cube complex} is built from gluing unit cubes along their faces by isometries. The \textit{link} of a vertex $v$, $lk(v)$, in a cube complex is a simplicial complex whose vertices lie on the edges adjacent to $v$. A collection of vertices in $lk(v)$ span a $(d - 1)$-dimensional complex if and only if they lie on edges of a common $d$-dimensional complex. A link $lk(v)$ is \textit{flag} if for every complete graph in $lk(v)$, there is a simplex in $lk(v)$ whose 1-skeleton is the complete graph. A cube complex is NPC if it satisfies a combinatorial link condition: all vertex links of the cube complex are flag complexes. We will later show that simple, non-positively curved surface amalgams are NPC cube complexes. 

In their seminal paper \cite{hw}, Haglund and Wise define \textit{special cube complexes}. We refer the reader to Section 3 of their paper for precise definitions, but special cube complexes satisfy three conditions: (1) every hyperplane embeds; (2) no hyperplane directly self-osculates; and (3) no two hyperplanes inter-osculate. A group $G$ is \textit{virtually special} if has a finite-index subgroup which is the fundamental group of a special cube complex. One main result in \cite{hw} is that Gromov hyperbolic virtually special groups are QCERF: 

\begin{theorem}[\cite{hw}, Corollary 7.4] \label{thm:hw} Let $X$ be a compact cube complex with Gromov hyperbolic fundamental group. If $X$ is virtually special, then every quasiconvex subgroup of $\pi_1(X)$ is separable (i.e., $\pi_1(X)$ is QCERF). 
\end{theorem}

Finally, we mention a convenient criterion developed by Wise for determining whether a Gromov hyperbolic group is virtually special: 

\begin{definition}[\cite{wise}, Definition 11.5] \label{def:QVH} Let $\mathcal{QVH}$ denote the smallest collection of hyperbolic groups closed under the following four operations:

\begin{enumerate}
\item $1 \in \mathcal{QVH}$;
\item If $G = A \ast_C B$, $A, B \in \mathcal{QVH}$, and $C$ is finitely generated and quasiconvex in $G$, then $G \in \mathcal{QVH}$; 
\item If $G = A \ast_C$, $A \in \mathcal{QVH}$, and $C$ is finitely-generated and quasiconvex in $G$, then $G \in \mathcal{QVH}$; 
\item Let $H \subset G$ be a finite-index subgroup and let $H \in \mathcal{QVH}$. Then $G \in \mathcal{QVH}$. 
\end{enumerate}
\end{definition} 

 The following deep theorem relates the hierarchy of groups $\mathcal{QVH}$ to virtual specialness:  

\begin{theorem}[\cite{wise}, Theorem 13.3] \label{thm:qvh} A torsion-free, Gromov hyperbolic group $G$ is virtually special if and only if $G \in \mathcal{QVH}$. 
\end{theorem}

Again, it is unknown which fundamental groups of surface amalgams are special; however, using \Cref{thm:qvh}, we will prove that fundamental groups of simple negatively curved surface amalgams are virtually special, and therefore QCERF. 
\section{The Base Case} 
\label{easy}

Recall that for any proper metric space $X$, the identity map $(X, g_1) \rightarrow (X, g_2)$ lifts to a quasi-isometry $f: (\widetilde{X}, \widetilde{g_1}) \rightarrow (\widetilde{X}, \widetilde{g_2})$ by the \u{S}varc-Milnor Lemma. By the Morse Lemma, $f$ in turn induces a boundary homeomorphism $\partial^{\infty}f: \partial^{\infty}(\widetilde{X}, \widetilde{g_1}) \rightarrow \partial^{\infty}(\widetilde{X}, \widetilde{g_2})$. The following section is devoted to proving $\partial^{\infty}f$ induces an isometry, a fact summarized in the following proposition: 

\begin{proposition}
\label{main1}
Suppose $X$ is a simple P-manifold endowed with a pair of metrics $g_1$ and $g_2 \in \mathcal{M}_{\leq}$. Furthermore, suppose that there is an inclusion of every chamber $C$ into a closed surface $S \subset X$. If $(X, g_1)$ and $(X, g_2)$ have the same marked length spectrum, then there exists an isometry $\phi: (X, g_1) \rightarrow (X, g_2)$ that is induced by the boundary homeomorphism $\partial^{\infty}f: \partial^{\infty}(\widetilde{X}, \widetilde{g_1}) \rightarrow \partial^{\infty}(\widetilde{X}, \widetilde{g_2})$ discussed above.
\end{proposition}

We first provide an outline of the proof. First, in the proof of Proposition \ref{mobius} of Section \ref{boundarymobius}, we show that the boundary homeomorphism $\partial^{\infty}f: \partial^{\infty}(\widetilde{X}, \widetilde{g_1}) \rightarrow \partial^{\infty}(\widetilde{X}, \widetilde{g_2})$ induced by the identity map is M\"{o}bius (see Definition \ref{crm}). The proof of Proposition \ref{mobius} requires the ergodicity of the geodesic flow map, which can be deduced from general theory developed by Kaimanovich in \cite{kai} for Gromov hyperbolic spaces (see the Appendix). Using ergodicity of the geodesic flow map, the proof of Lemma \ref{approx} shows the cross ratio of any distinct $(a, b, c, d) \in \big(\partial^{\infty}(\widetilde{X})\big)^4$ can be approximated arbitrarily well by lengths of closed geodesics in $(X, g)$. We can thus conclude that $\partial^{\infty}f$ is M\"{o}bius since $(X, g_1)$ and $(X, g_2)$ have the same marked length spectra, which determine lengths of closed geodesics and thus cross ratios. Next, in Section \ref{patching}, we patch together isometries $\phi_S: (S, g_1\vert_{S}) \rightarrow (S, g_2\vert_{S})$ constructed by Otal in \cite{Otal90} on every closed subsurface $S \subset X$ to construct a global isometry $\phi: (X, g_1) \rightarrow (X, g_2)$. In particular, we show that given two surfaces $S$ and $S'$ that intersect in $X$, $\phi_S\vert_{S \cap S'} = \phi_{S'}\vert_{S \cap S'}$ using the fact that $\partial^{\infty}f$ is M\"{o}bius. 

\subsection{The boundary map is M\"{o}bius}
\label{boundarymobius}

In this section, we prove the following proposition by adapting Otal's proof of Theorem 1 in \cite{Otal90}. Note that in this section, the metrics are only assumed to be locally CAT(-1). 

\begin{proposition}\label{mobius} Suppose $(X, g_1)$ and $(X, g_2)$ are two simple, thick locally CAT(-1) P-manifolds. If $(X, g_1)$ and $(X, g_2)$ have the same marked length spectrum, then the boundary map is M\"{o}bius. 
\end{proposition}

A key ingredient in Otal's proof is ergodicity of the geodesic flow map on surfaces, a well-known and classical result. Recall that in the setting of surfaces, geodesic flow is defined on the unit tangent bundle. For P-manifolds, however, the unit tangent bundle is undefined on the gluing curves. One, however, has the notion of a widely-studied \textit{generalized unit tangent bundle}, which can be identified with the usual unit tangent bundle in the setting of Riemannian manifolds. 

Recall that for a metric space $(X, g)$, a geodesic $\gamma: \mathbb{R} \rightarrow X$ has \textit{speed} $s \geq 0$ if for every $t \in \mathbb{R}$, there exists a neighborhood $U \subset \mathbb{R}$ such that for all $t_1, t_2 \in U$, $g\big(\gamma(t_1), \gamma(t_2)\big) = s\lvert t_1 - t_2 \rvert$. In particular, if $s = 1$, $\gamma$ is a \textit{unit-speed geodesic}. 

\begin{definition}[Generalized unit tangent bundle] Given a geodesically complete metric space $X$, the \textit{generalized unit tangent bundle} of $X$, $SX$, is the space of unit-speed geodesics in $X$. 
\end{definition} 

Suppose two unit-speed geodesics $\gamma \sim \gamma'$ if $\gamma(t) = \gamma'(-t)$ for all $t \in \mathbb{R}$. Then there is a natural identification of $SX/\sim$ with $\mathscr{G}(X) \times \mathbb{R}$, where as before, $\mathscr{G}(X)$ denotes the space of unoriented, unparametrized geodesics in $X$. We say a point $x \in X$ is a \textit{basepoint} of $\xi \in SX$ if $x = \xi(0)$. We now recall the definition of the geodesic flow map on $SX$:  

\begin{definition}[Geodesic flow] The \textit{geodesic flow} on $SX$ is the map $\phi_t: SX \rightarrow SX, \phi_t(\xi)(s) = \xi(t + s)$, where $s, t \in \mathbb{R}$. 
\end{definition}

For a closed negatively curved Riemannian surface, there is a unique direction in which to continue a geodesic via the exponential map. On the other hand, for a P-manifold, there are $n - 1$ directions in which to continue a geodesic hitting a point on a gluing curve $\gamma$, where $n$ is the number of boundary components glued to $\gamma$. Due to this key difference between surfaces and P-manifolds, the usual proof of the ergodicity of the geodesic flow map on surfaces does not generalize to the setting of P-manifolds and different techniques are needed. We refer the reader to the Appendix, which summarizes a general result by Kaimanovich about the ergodicity of the geodesic flow map on proper Gromov hyperbolic spaces. Kaimanovich's results apply in our setting since simple, thick P-manifolds equipped with locally CAT(-1) metrics are proper Gromov hyperbolic spaces. 

Recall from \cite{bk2} that a sequence of points $(x_n)_{n \in \mathbb{N}}$ in a Gromov hyperbolic space $X$ \textit{converges to infinity} if for any choice of $x \in X$, $$\lim \inf\limits_{i, j \rightarrow \infty} \langle x_i, x_j \rangle_{x} = \infty$$ where $\langle \cdot \rangle_{x}$ denotes the Gromov Product from Definition \ref{gromovproduct}; this definition is independent of choice of basepoint $x$. Furthermore, two sequences $(x_n)_{n \in \mathbb{N}}$ and $(y_n)_{n \in \mathbb{N}}$ are equivalent if
$$\lim \inf\limits_{i, j \rightarrow \infty} \langle x_i, y_j \rangle_{x} = \infty.$$ Endowing $\overline{X} = X \cup \partial^{\infty}(X)$ with the topology described in Definition 2.13 of \cite{bk2}, we can think of points in $\partial^{\infty}$ as limit points of equivalence classes of sequences of points in $X$ that converge to infinity.

With the above paragraph in mind, we introduce a lemma whose proof is adapted from Otal's proof of Theorem 1 in \cite{otal92}:

\begin{lemma} \label{approx} Suppose $X$ is simple, thick P-manifold endowed with a locally CAT(-1) metric. Then given a 4-tuple of distinct points $(a, b, c, d) \in (\partial^{\infty}(\widetilde{X}))^4$, $[abcd]$ can be approximated arbitrarily well by lengths of closed geodesics.   
\end{lemma}

\begin{proof} Fix a 4-tuple of distinct points $(a, b, c, d) \in \big(\partial^{\infty}(\widetilde{X}, \widetilde{g})\big)^4$. By the Birkhoff Ergodic Theorem (see Section \ref{birkhoffapp}), we can conclude that there exists $v \in SX$ with dense orbit under the geodesic flow map $\phi_t$ on $SX$. Thus, there exist sequences of vectors $(\phi_{n}(\widetilde{v_1}))_{n \in \mathbb{N}}, (\phi_{-n}(\widetilde{v_1}))_{n \in \mathbb{N}}, (\phi_{n}(\widetilde{v_2}))_{n \in \mathbb{N}},$ and $(\phi_{-n}(\widetilde{v_2}))_{n \in \mathbb{N}}$ whose basepoints approach $a, b, c,$ and $d$ respectively (here, $\widetilde{v_1}$ and $\widetilde{v_2}$ are lifts of $v$ in $S\widetilde{X}$). As a result, Then there exist $n_i \in \mathbb{N}$ such that $[\phi_{n_1}(\widetilde{v_1})(0)\phi_{-n_2}(\widetilde{v_1})(0)\phi_{n_3}(\widetilde{v_2})(0)\phi_{-n_4}(\widetilde{v_2})(0)]$ approximates $[abcd]$ arbitrarily well. 


Fix any $u \in S\widetilde{X}$. Then since the orbit of $v$ is dense in $SX = S\widetilde{X} / \Gamma$, for every sufficiently large $m_1 \in \mathbb{N}$ there exists some $\gamma_1 \in \Gamma$ such that $\gamma_1\cdot\phi_{m_1}(\widetilde{v_1})$ lies arbitrarily close to $\gamma_1 \cdot u$ in $S\widetilde{X} / \Gamma$. Similarly, there exists some $\gamma_3 \in \Gamma$ and sufficiently large $m_3 \in \mathbb{N}$ such that $\gamma_3 \cdot \phi_{m_3}(\widetilde{v_2})$ and $\gamma_3 \cdot u$ are arbitrarily close. From this, by choosing $\gamma_l, \gamma'_m \in \Gamma$ accordingly, we can construct sequences of vectors $(\gamma_m \cdot u)_{m \in \mathbb{N}}$ and $(\gamma'_l \cdot u)_{l \in \mathbb{N}}$ whose basepoints limit to $a$ and $c$ respectively (see Figure \ref{fig:approx2}). Note that the translation lengths of $\gamma_l$ and $\gamma'_m$ will become arbitrarily large.


 \begin{figure}[h!]
 \begin{centering}
\includegraphics[width=0.4\textwidth]{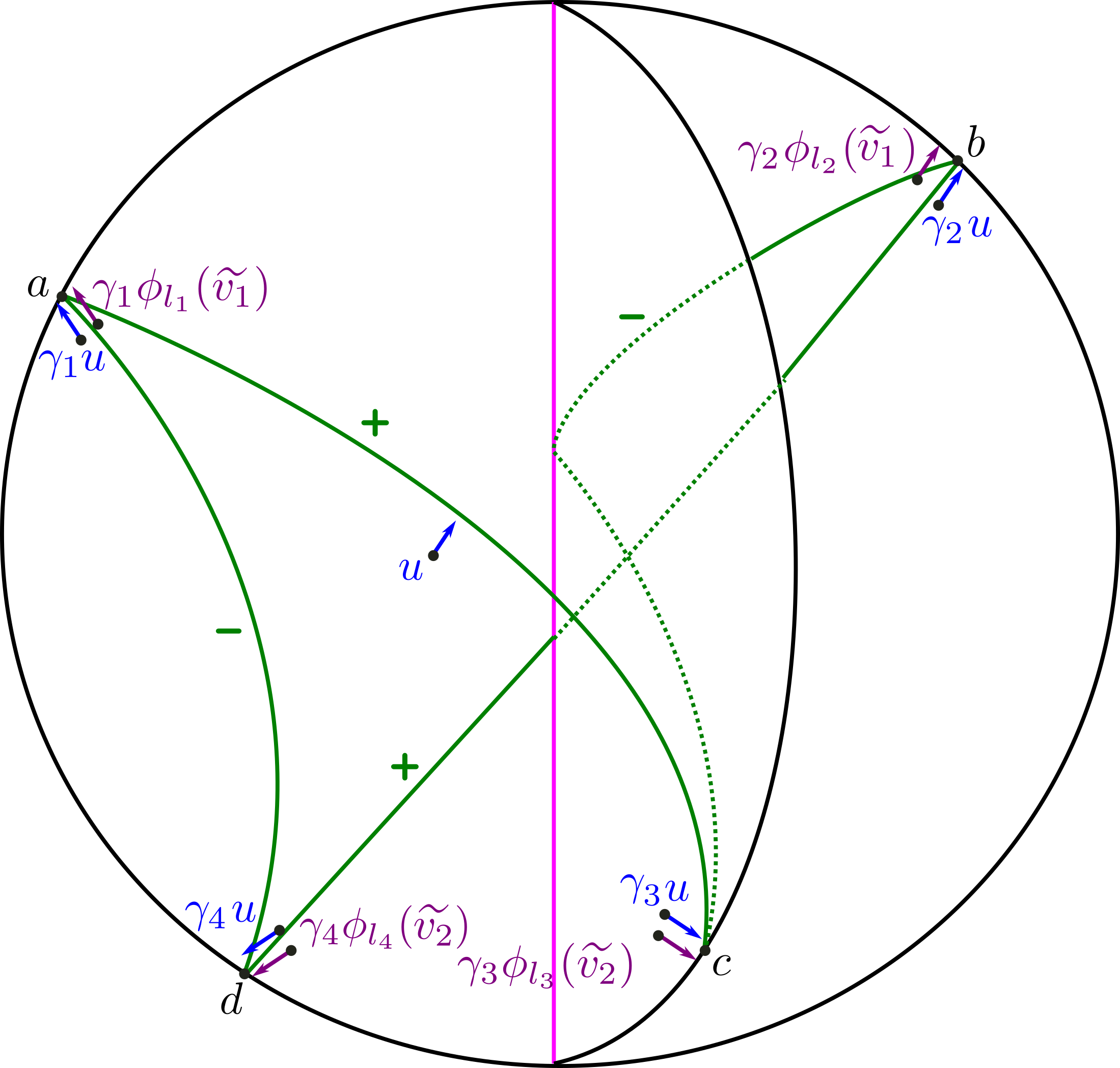}
    \caption{An illustration of the setup of Lemma \ref{approx}.}
    \label{fig:approx2}
\end{centering} 
\end{figure}
 
 Let $\ell_{\gamma^{-1}_m\gamma'_l}$ be the length of the fundamental domain of the action of $\gamma^{-1}_m \circ \gamma'_l$ on its translation axis. Note that $\ell_{\gamma^{-1}_m \gamma'_l}$, the length of a closed geodesic corresponding to $\gamma^{-1}_m \circ \gamma'_l$, is approximated arbitrarily well by $\widetilde{g}\big((\gamma_m \cdot u)(0), (\gamma'_l \cdot u)(0)\big)$ since by construction, the basepoints of $(\gamma_m \cdot u)_{m \in \mathbb{N}}$ and $(\gamma'_l \cdot u)_{l \in \mathbb{N}}$ approach the endpoints of a geodesic $(a, c) \in (\widetilde{X}, \widetilde{g})$ and $\gamma_1 \cdot u = (\gamma_1 \circ \gamma_3^{-1})(\gamma_3 \cdot u)$. We apply a similar argument to construct sequences $(\gamma''_k \cdot u)_{k \in \mathbb{N}}$ and $(\gamma'''_j \cdot u)_{j \in \mathbb{N}}$ that tend to $b$ and $d$ respectively. Furthermore, we can construct the sequences such that $\widetilde{g}$-distances between basepoints of terms of the sequences are arbitrarily close to lengths of long closed geodesics. 
 
 Since $(\gamma_m \cdot u)_{m \in \mathbb{N}}$, $(\gamma''_k \cdot u)_{k \in \mathbb{N}}$, $(\gamma'_l \cdot u)_{l \in \mathbb{N}}$, and $(\gamma'''_j \cdot u)_{j \in \mathbb{N}}$ are sequences whose basepoints tend to $a$, $b$, $c$, and $d$ respectively, we have that $[abcd]$ is approximated arbitrarily well by $[(\gamma_m \cdot u)(0) (\gamma''_k \cdot u)(0) (\gamma'_l \cdot u)(0) (\gamma'''_j \cdot u)(0)]$. Since we can approximate $[(\gamma_m \cdot u)(0) (\gamma''_k \cdot u)(0) (\gamma'_l \cdot u)(0) (\gamma'''_j \cdot u)(0)]$ by adding and subtracting lengths of closed geodesics, the result then follows. 
 

\end{proof}

We can now prove Proposition \ref{mobius}. By Lemma \ref{approx}, for all $i \in \mathbb{N}$, there is some quadruple of geodesic arcs $\alpha_i = (x_i, y_i)$, $\beta_i = (x'_i, y'_i)$, $\gamma_i = (x_i, y'_i)$, and $\delta_i = (x'_i, y_i)$ projecting to closed geodesics in $X$ such that:

$$\big\lvert [abcd] - \big(\widetilde{g_1}(x_i, y_i) + \widetilde{g_1}(x'_i, y'_i) - \widetilde{g_1}(x_i, y'_i) - \widetilde{g_1}(x'_i, y_i)\big) \big\rvert < \frac{1}{i}.$$

By construction, $(x_i, y_i, x'_i, y'_i)$ converges to $(a, b, c, d)$ in $(\widetilde{X}, \widetilde{g_1})$ so $(f(x_i), f(y_i), f(x'_i), f(y'_i))$ will converge to $(\partial^{\infty}f(a), \partial^{\infty}f(b), \partial^{\infty}f(c), \partial^{\infty}f(d))$ in $(\widetilde{X}, \widetilde{g_2})$ since $f$ is continuous. Thus, by definition of cross ratio, we have that:
$$[abcd] = \lim\limits_{i \rightarrow \infty} \widetilde{g_1}(x_i, y_i) + \widetilde{g_1}(x'_i, y'_i) - \widetilde{g_1}(x_i, y'_i) - \widetilde{g_1}(x'_i, y_i)$$
and
$$\lim\limits_{i \rightarrow \infty} \widetilde{g_2}(f(x_i, y_i)) + \widetilde{g_2}(f(x'_i, y'_i)) - \widetilde{g_2}(f(x_i, y'_i)) - \widetilde{g_2}(f(x'_i, y_i)) = [\partial^{\infty}f(a) \partial^{\infty}f(b) \partial^{\infty}f(c) \partial^{\infty}f(d)].$$
Recall that since $(x_i, y_i)$, $(x'_i, y'_i)$, $(x'_i, y_i)$ and $(x_i, y'_i)$ each project to closed geodesics in $X$, the distances between their endpoints are all lengths of elements in $\pi_1(X)$. Thus, since $(X, g_1)$ and $(X, g_2)$ have the same marked length spectrum, $\widetilde{g_1}(x_i, y_i) = \ell_{g_1}(\alpha_i) = \ell_{g_2}(\alpha_i) = \widetilde{g_2}(f(x_i, y_i))$, and a similar statement holds for the other distances as well. It then follows that $[abcd] = [\partial^{\infty}f(a) \partial^{\infty}f(b) \partial^{\infty}f(c) \partial^{\infty}f(d)]$. \qed 

\subsection{Proof of Proposition \ref{main1}}
\label{patching} 

In order to prove Proposition \ref{main1}, we need a few auxiliary lemmas. The first lemma is a basic fact about patching isometries together. 

Recall a metric space $(X, d)$ is \textit{convex} if for any two points $x, y \in X$, there exists $z \in X$ distinct from $x$ and $y$ such that $d(x, z) + d(z, y) = d(x, y)$. 

\begin{lemma}\label{isom}
Suppose $U_i$ and $V_i$ are complete, convex, locally compact metric spaces, where $i = 1, 2$. Suppose $\phi_i: U_i \rightarrow V_i$ are (invertible) isometries, and ${\phi_1}\vert_{U_1 \cap U_2} = {\phi_2}\vert_{U_1 \cap U_2}$. Then there exists an isometry $\phi: U_1 \cup U_2 \rightarrow V_1 \cup V_2$ such that $\phi\vert_{U_i} = {\phi_i}\vert_{U_i}$ for $i = 1, 2$. 
\end{lemma}

\begin{proof} In the following proof, for a metric space $X$, $d_X$ will denote the distance function of $X$. We define $\phi$ in the natural way: $\phi(x) = \phi_i(x)$ if $x \in U_i$. It suffices to show the theorem is true when $p$ and $q$ are not both in $U_1$ or $U_2$. Suppose without loss of generality $p \in U_1$ and $q \in U_2$. By Hopf-Rinow, since $U_1 \cup U_2$ is complete, convex, and locally compact, then there exists a minimizing geodesic between $p$ and $q$, which we will call $[p, q]$. Let $r \in (U_1 \cap U_2) \cap [p, q]$, so $d_{U_1 \cup U_2}(p, r) + d_{U_1 \cup U_2}(r, q) = d_{U_1 \cup U_2}(p, q)$ (see Figure \ref{fig:isom}).

\begin{figure}[h!]
    \centering
    \includegraphics[width=0.9\textwidth]{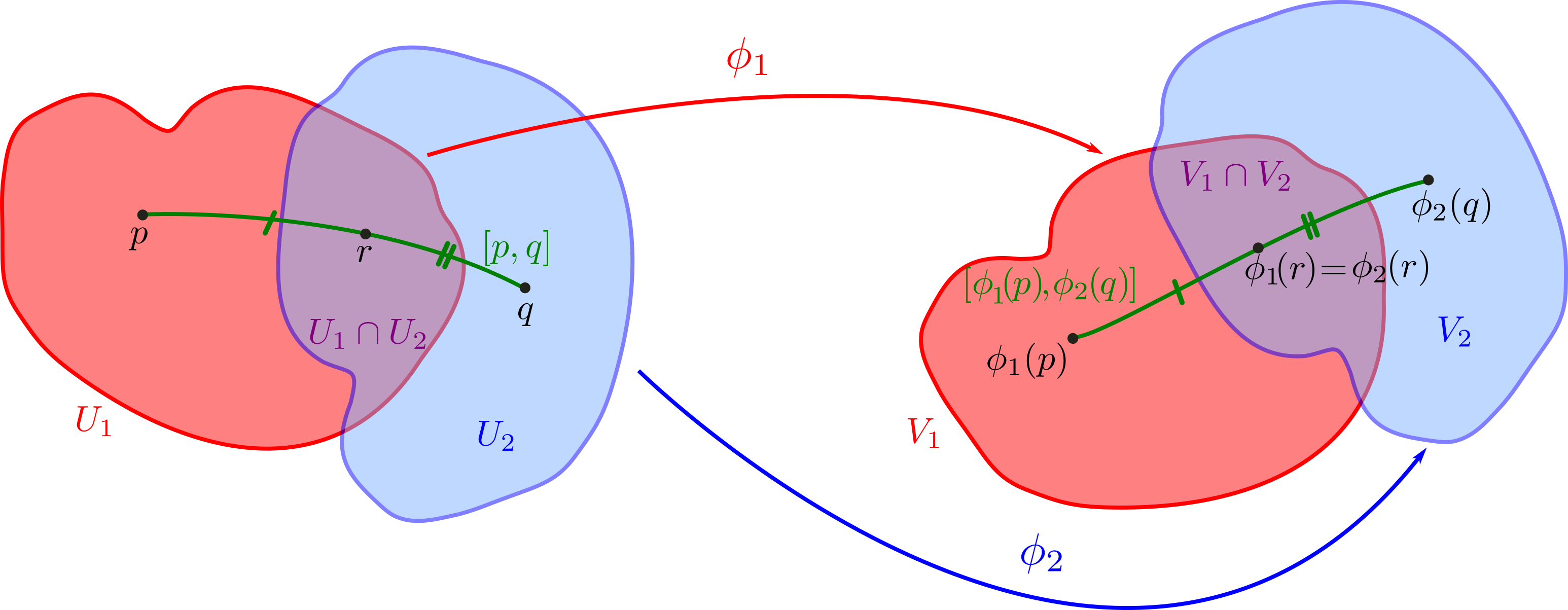}
    \caption{An illustration of the setup of Lemma \ref{isom}.}
    \label{fig:isom}
\end{figure}

Since $\phi_i$ is an isometry for $i = 1, 2$, it follows that: 
\begin{align*}
d_{V_1 \cup V_2}(\phi(p), \phi(r)) + d_{V_1 \cup V_2}(\phi(r), \phi(q)) &= d_{V_1}(\phi_1(p), \phi_1(r)) + d_{V_2}(\phi_2(r), \phi_2(q))\\
&= d_{U_1}(p, r) + d_{U_2}(r, q) = d_{U_1 \cup U_2}(p, q).
\end{align*}
By triangle inequality, $d_{V_1 \cup V_2}(\phi(p), \phi(q)) \leq d_{U_1 \cup U_2}(p, q)$. For the reverse direction, apply a symmetric argument: since $V_1 \cup V_2$ is complete, convex, and locally compact, there exists a minimizing geodesic between $\phi(p)$ and $\phi(q)$, $[\phi(p), \phi(q)]$. Choose $r' \in [\phi(p), \phi(q)] \cap (V_1 \cap V_2)$ such that $d_{V_1}(\phi(p), r') + d_{V_2}(r', \phi(q)) = d_{V_1 \cup V_2}(\phi(p), \phi(q))$. Since $\phi_i$ is invertible for $i = 1, 2$, 
\begin{align*}d_{V_1 \cup V_2}(\phi(p), \phi(q)) &=  d_{V_1 \cup V_2}(\phi_1(p), \phi_2(q)) = d_{V_1}(\phi_1(p), r') + d_{V_2}(r', \phi_2(q)) \\
&= d_{U_1}(p, \phi_1^{-1}(r')) + d_{U_2}(\phi_2^{-1}(r'), q) = d_{U_1 \cup U_2}(p, \phi^{-1}(r')) + d_{U_1 \cup U_2}(\phi^{-1}(r'), q)\\
&\geq d_{U_1 \cup U_2}(p, q).
\end{align*}

Then it follows $d_{V_1 \cup V_2}(\phi(p),\phi(q)) = d_{U_1 \cup U_2}(p, q)$.  
\end{proof}

The following is a technical lemma about the cross-ratio defined in Equation \ref{cr}.  

\begin{lemma}
\label{crossratio} 
Suppose $(\widetilde{S}, \widetilde{g}\vert_{\widetilde{S}})$ and $(\widetilde{S'}, \widetilde{g}\vert_{\widetilde{S'}})$ are two lifts of $(S, g\vert_S)$ and $(S', g\vert_{S'})$ respectively that meet at a branching geodesic $\widetilde{\gamma} \subset (\widetilde{X}, \widetilde{g})$. Then for every $p \in \widetilde{\gamma}$, there exist four points $a, b, c, d$ such that the geodesics $(a, c)$, $(b, d)$, $(a, d)$, and $(b, c)$ all meet at $p$ and $[abcd] = 0$.  
\end{lemma}

\begin{proof} Choose some $p \in \widetilde{\gamma}$. We will find a set of four geodesics $(a, c)$, $(b, d)$, $(a, d)$, and $(b, c)$ that all meet at $p$. Choose an arbitrary $(a, c) \subset (\widetilde{S}, \widetilde{g}\vert_{\widetilde{S}})$ that transversely intersects $\widetilde{\gamma}$ at $p$ in $\widetilde{S}$. Suppose $\widetilde{\gamma}$ divides $\widetilde{S'}$ into two half planes, which we will call $H'_1$ and $H'_2$. Consider a map $h: \partial^{\infty}(H'_1) \rightarrow \widetilde{\gamma}$ that maps a point $x \not\in \partial^{\infty}(H'_1)$ to the point $h(x) = \widetilde{\gamma} \cap (x, c)$ where $(x, c)$ denotes the geodesic arc with endpoints $x$ and $c$. Note that $h$ extends continuously to a map $\widetilde{h}: \partial^{\infty}(H'_1) \rightarrow \widetilde{\gamma} \cup \{q, r\}$, where $q, r \in \partial^{\infty}(\widetilde{S}) \cap \partial^{\infty}(\widetilde{S'})$ are endpoints of $\widetilde{\gamma}$ in the following way: $\widetilde{h}(q) = q$ and $\widetilde{h}(r) = r$. As we vary $x$ along $\partial^{\infty}(H'_1)$, $h(x)$ varies continuously along $\widetilde{\gamma}$, so by an Intermediate Value Theorem argument, since $\widetilde{\gamma}$ and $\partial^{\infty}(H'_1)$ are both connected, it follows that $\widetilde{h}(x) = p$ for some $x \in \partial^{\infty}(H'_1)$ (see Figure \ref{fig:ivt}). We set $b$ equal to this $x$. 

\begin{figure}[h!]
    \centering
    \includegraphics[width=0.4\textwidth]{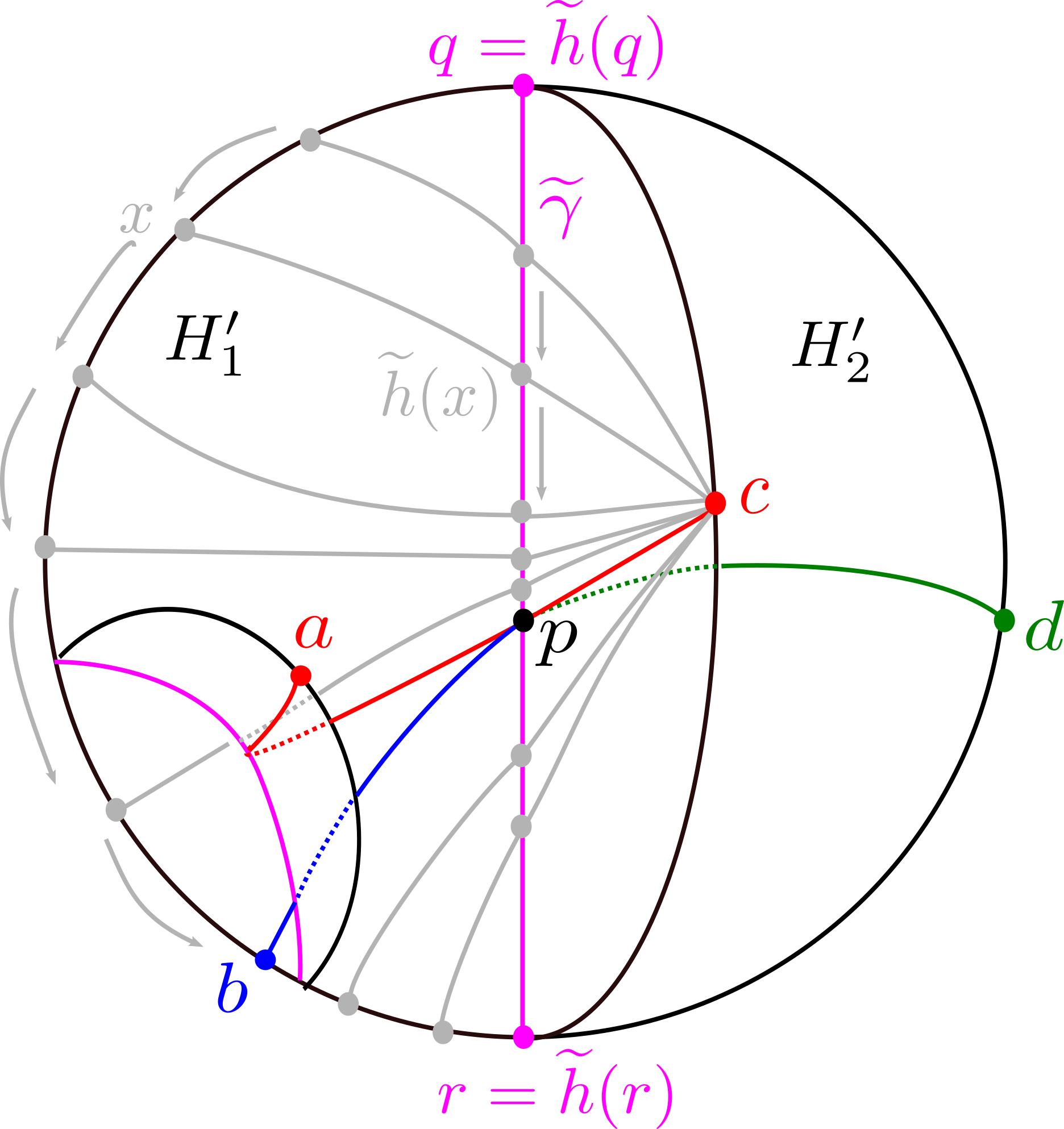}
    \caption{By an intermediate value theorem argument, there exists some $x \in \partial^{\infty}H'_1$ such that $\widetilde{h}(x) = p$ for a given point $p \in \widetilde{\gamma}$.}
    \label{fig:ivt}
\end{figure}

Note that $(b, c)$ is a geodesic that passes through $p$. Continue the geodesic ray $(b, p)$ into $H'_2 \subset \widetilde{S'}$ to obtain a bi-infinite geodesic $(b, d) \subset \widetilde{S'
}$. Note that $(a, d)$ will also be a geodesic that passes through $p$. Thus, we have found four geodesics $(a, c)$, $(b, d)$, $(a, d)$, and $(b, c)$ that all intersect at $p$. It is then possible to find a sequence $(a_i, b_i, c_i, d_i)$ where $a_i \in (p, a)$, $b_i \in (p, b)$, $c_i \in (p, c)$ and $d_i \in (p, d)$ such that: 

\begin{dmath*}
[abcd] = \lim\limits_{(a_i, b_i, c_i, d_i) \rightarrow (a, b, c, d)}\widetilde{g_1}(a_i, c_i) + \widetilde{g_1}(b_i, d_i) - \widetilde{g_1}(a_i, d_i) - \widetilde{g_1}(b_i, c_i) 
= \lim\limits_{(a_i, b_i, c_i, d_i) \rightarrow (a, b, c, d)} \widetilde{g_1}(a_i, p) + \widetilde{g_1}(p, c_i) + \widetilde{g_1}(b_i, p) + \widetilde{g_1}(p, d_i) - \widetilde{g_1}(a_i, p) - \widetilde{g_1}(p, d_i) - \widetilde{g_1}(b_i, p) - \widetilde{g_1}(p, c_i)\\
= 0. 
\end{dmath*}

\end{proof}

Let $\mathcal{S}$ be a collection of closed surfaces that covers $X$. We know such a covering exists; for each chamber $C_i \subset X$, note that $C_i \subset S_i$ by assumption, where $S_i$ is a closed surface. Then $\{S_i\}_{i = 1}^{n}$ is a collection of closed surfaces that covers $X$.  

\begin{lemma}
\label{chamberisom}
Suppose $S, S' \in \mathcal{S}$ are two closed surfaces in $X$ that are identified along some set of gluing curves $\{\gamma_1, ..., \gamma_n\}$. Then there exist isometries $\phi_{S}: (S, g_1\vert_{S}) \rightarrow (S, g_2\vert_{S})$ and $\phi_{S'}: (S', g_1\vert_{S'}) \rightarrow (S', g_2\vert_{S'})$ where $\phi_{S}\vert_{\gamma_i} = \phi_{S'}\vert_{\gamma_i}$ for all $1 \leq i \leq n$. 
\end{lemma}

\begin{proof}
Suppose $(X, g_1)$ and $(X, g_2)$ have the same marked length spectra. As before, $\partial^{\infty}f$ will denote the boundary homeomorphism between $\partial^{\infty}(\widetilde{X}, \widetilde{g_1})$ and $\partial^{\infty}(\widetilde{X}, \widetilde{g_2})$. 

Recall $\widetilde{S}$ and $\widetilde{S'}$ are apartments that are arbitrary lifts of $S$ and $S'$ in $\widetilde{X}$ such that $\widetilde{S} \cap \widetilde{S'} \neq \varnothing$. We can then use the restrictions of the boundary homeomorphism $\partial^{\infty}f\vert_{\partial^{\infty}(\widetilde{S})}$ and $\partial^{\infty}f\vert_{\partial^{\infty}(\widetilde{S}')}$ to construct isometries $\phi_{S}: (S, g_1\vert_{S}) \rightarrow (S, g_2\vert_{S})$ and $\phi_{S'}: (S', g_1\vert_{S'}) \rightarrow (S', g_2\vert_{S'})$ via the methods of Otal. In other words, given $x \in (\widetilde{S}, \widetilde{g_1}\vert_{\widetilde{S}})$, consider the set of all geodesics in $(\widetilde{S}, \widetilde{g_1}\vert_{\widetilde{S}})$ that intersect at $x$. By \cite{Otal90}, the set of geodesics mapped to $(\widetilde{S}, \widetilde{g_2}\vert_{\widetilde{S}})$ via the boundary map $\partial^{\infty}f\vert_{\partial^{\infty}(\widetilde{S})}$ intersect at a single point $\widetilde{\phi_{S}}(x)$; $\widetilde{\phi}_{S'}$ is constructed similarly. Otal shows that $\widetilde{\phi_S}$ and $\widetilde{\phi_{S'}}$ are $\pi_1$-equivariant. 

We first check that $\widetilde{\phi_{S}}$ and $\widetilde{\phi_{S'}}$ extend to $\partial^{\infty}f\vert_{\partial^{\infty}(\widetilde{S})}$ and $\partial^{\infty}f\vert_{\partial^{\infty}(\widetilde{S'})}$ respectively. Indeed, consider a sequence of points $\{x_i\}_{i \in \mathbb{N}}$ that converge radially towards a point $x \in \partial^{\infty}\widetilde{S}$, so $x = \lim\limits_{i \rightarrow \infty} x_i$. It is possible to find a sequence of pairs of geodesics $(\xi_i, \eta_i)$, such that $\xi_i \cap \eta_i = x_i$ and $a_i, b_i, c_i, d_i \rightarrow x$, where $a_i$ and $b_i$ are the endpoints of $\xi_i$ and $c_i$ and $d_i$ are the endpoints of $\eta_i$ (see Figure \ref{fig:radial}). Then since $\partial^{\infty}f$ is continuous, $\lim\limits_{i \rightarrow \infty} \partial^{\infty}f(a_i) = \partial^{\infty}f(x)$, and the same is true for $b_i, c_i$, and $d_i$. Then: 
$$\lim\limits_{i \rightarrow \infty}\widetilde{\phi_S}(x_i) = \lim\limits_{i \rightarrow \infty} \widetilde{\phi_S}\big((a_i, b_i) \cap (c_i, d_i)\big) = \lim\limits_{i \rightarrow \infty} \big( \partial^{\infty}f(a_i), \partial^{\infty}f(b_i)\big) \cap \big( \partial^{\infty}f(c_i), \partial^{\infty}f(d_i)\big) =  \partial^{\infty}f(x).$$

\begin{figure}[h!]
    \centering
    \includegraphics[width=0.35\textwidth]{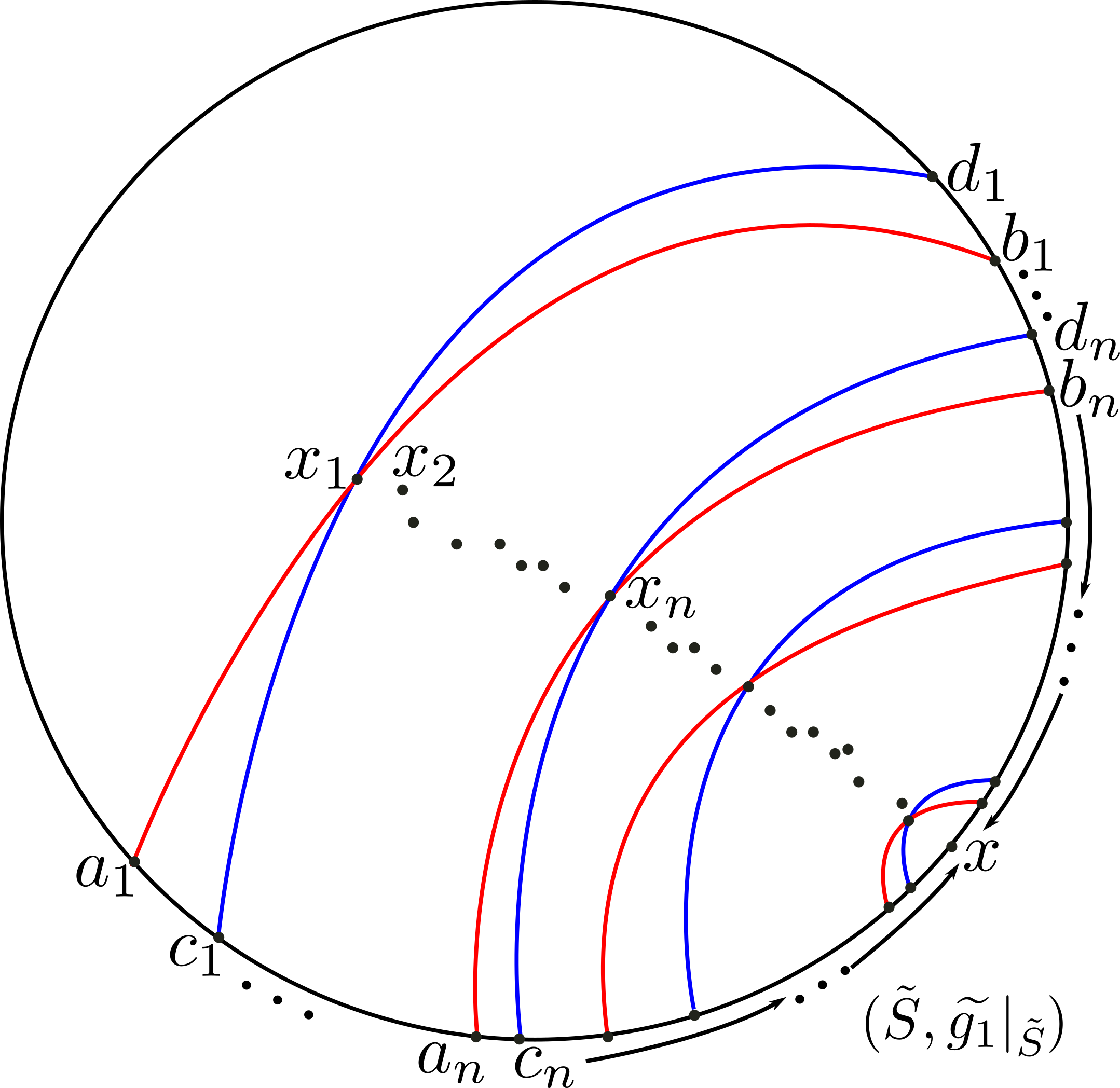}
    \caption{Given a sequence of points $\{x_i\}_{i = 1}^n$ converging radially to $x$, one can find a sequence of geodesics $\xi_i$ and $\eta_i$ with endpoints $a_i, b_i, c_i, d_i$ converging to $x$. The images of these endpoints under the boundary homeomorphism $\partial^{\infty}f$ will also converge to $\partial^{\infty}f(x) = 
    \lim\limits_{i \rightarrow \infty} \widetilde{\phi_S}(x_i)$. }
    \label{fig:radial}
\end{figure}

As before, $\big(\partial^{\infty}f(a_i), \partial^{\infty}f(b_i)\big)$ denotes a geodesic with endpoints $\partial^{\infty}f(a_i)$ and $\partial^{\infty}f(b_i)$. 

As a result, we conclude that $\partial^{\infty}f\vert_{\partial^{\infty}(\widetilde{S})} = \partial^{\infty}\widetilde{\phi_S}$ , $\partial^{\infty}f\vert_{\partial^{\infty}(\widetilde{S'})} = \partial^{\infty}\widetilde{\phi_{S'}}$, and $\partial^{\infty}\widetilde{\phi_S}\vert_{\widetilde{S} \cap \widetilde{S'}} = \partial^{\infty}f\vert_{\partial^{\infty}(\widetilde{S} \cap \widetilde{S'})} = \partial^{\infty}\widetilde{\phi_{S'}}\vert_{\widetilde{S} \cap \widetilde{S'}}$. In particular, if $p \in \partial^{\infty}(\widetilde{S}) \cap \partial^{\infty}(\widetilde{S'})$, then: 
\begin{equation}
\label{endpoints}
\partial^{\infty}\widetilde{\phi_S}(p) = \partial^{\infty}f\vert_{\widetilde{S} \cap \widetilde{S'}}(p) = \partial^{\infty}\widetilde{\phi_{S'}}(p). 
\end{equation}

Since isometries between Riemannian manifolds necessarily preserve geodesics, $\widetilde{\phi_S}$ and $\widetilde{\phi_{S'}}$ must map the branching geodesic $\widetilde{\gamma} = (p, q) \subset (\widetilde{S} \cap \widetilde{S'}, \widetilde{g_1})$ to geodesics in $(\widetilde{S}, \widetilde{g_2}\vert_{\widetilde{S}})$ and $(\widetilde{S'}, \widetilde{g_2}\vert_{\widetilde{S'}})$ respectively. By Equation \ref{endpoints}, $\widetilde{\phi_S}(\widetilde{\gamma})$ and $\widetilde{\phi_{S'}}(\widetilde{\gamma})$ must share endpoints. Since $\widetilde{S}, \widetilde{S'} \subset (\widetilde{X}, \widetilde{g_2})$ and geodesics between two given boundary points in the CAT(-1) space $(\widetilde{X}, \widetilde{g_2})$ are unique, $\widetilde{\phi_S}(\widetilde{\gamma}) = \widetilde{\phi_{S'}}(\widetilde{\gamma})$ necessarily. Thus, for every branching geodesic $\widetilde{\gamma} \subset (\widetilde{S}, \widetilde{g_1}) \cap (\widetilde{S'}, \widetilde{g_1})$, $\widetilde{\phi_S}(\widetilde{\gamma}) = \widetilde{\phi_{S'}}(\widetilde{\gamma})$.

It then suffices to show that for $x \in \widetilde{\gamma}$, $\widetilde{\phi_S}(x) = \widetilde{\phi_{S'}}(x)$. Suppose there exists some $x \in \widetilde{\gamma}$ where $\widetilde{\phi_S}(x) \neq \widetilde{\phi_{S'}}(x)$.
By Proposition \ref{mobius}, if an isomorphism of fundamental groups preserves the marked length spectrum of two CAT(-1) spaces, then the induced map at infinity is Mobius. As a consequence, $\partial^{\infty}f$ is Mobius. By Lemma \ref{crossratio}, for any $x \in \widetilde{\gamma}$ there exists a pair of bi-infinite geodesics $\xi = (a, b) \in (\partial^{\infty} \widetilde{S} \times \partial^{\infty} \widetilde{S}) \setminus \Delta$ and $\xi' = (a', b') \in (\partial^{\infty} \widetilde{S'} \times \partial^{\infty} \widetilde{S'}) \setminus \Delta$ such that $[aa'bb'] = 0$ (see Figure \ref{fig:branchingisom}). Note that $(a, b)$ and $(a', b')$ map to a pair of non-intersecting geodesics since we assumed that $\widetilde{\phi_S}(x) \neq \widetilde{\phi_{S'}}(x)$. 

Let $a_i, b_i \in \xi$, $a'
_i, b'_i \in \xi'$, $(a_i, b_i) \rightarrow (a, b)$ and $(a'_i, b'_i) \rightarrow (a', b')$. The Parallelogram Law for CAT($\kappa$) spaces (see Exercise 1.16 of \cite{bh}) gives us the following inequality: 
\begin{equation}
\label{parallelogram}
\widetilde{g_2}(\widetilde{\phi_S}(a'_i), \widetilde{\phi_S}(b_i)) + \widetilde{g_2}(\widetilde{\phi_S}(a_i), \widetilde{\phi_S}(b'_i)) \neq \widetilde{g_2}(\widetilde{\phi_S}(a_i), \widetilde{\phi_S}(b_i)) + \widetilde{g_2}(\widetilde{\phi_{S'}}(a'_i), \widetilde{\phi_{S'}}(b'_i)).
\end{equation} 

We therefore conclude that: 

\begin{align*}
[\partial^{\infty}f(a)\partial^{\infty}f(a')\partial^{\infty}f(b)\partial^{\infty}f(b')] &= \lim\limits_{i \rightarrow \infty}[\widetilde{\phi_S}(a_i) \widetilde{\phi_{S'}}(a'_i) \widetilde{\phi_S}(b_i) \widetilde{\phi_{S'}}(b'_i)]\\
& = \lim\limits_{i \rightarrow \infty} \widetilde{g_2}(\widetilde{\phi_S}(a_i), \widetilde{\phi_S}(b_i)) + \widetilde{g_2}(\widetilde{\phi_{S'}}(a'_i), \widetilde{\phi_{S'}}(b'_i)) \\
& - \widetilde{g_2}(\widetilde{\phi_S}(a_i), \widetilde{\phi_{S'}}(b'_i)) - \widetilde{g_2}(\widetilde{\phi_S}(b_i), \widetilde{\phi_{S'}}(a'_i)) \\
&\overset{(\ref{parallelogram})}{\neq} 0 = [aa'bb'].
\end{align*}

\begin{figure}[h!]
    \centering
\includegraphics[width=0.8\textwidth]{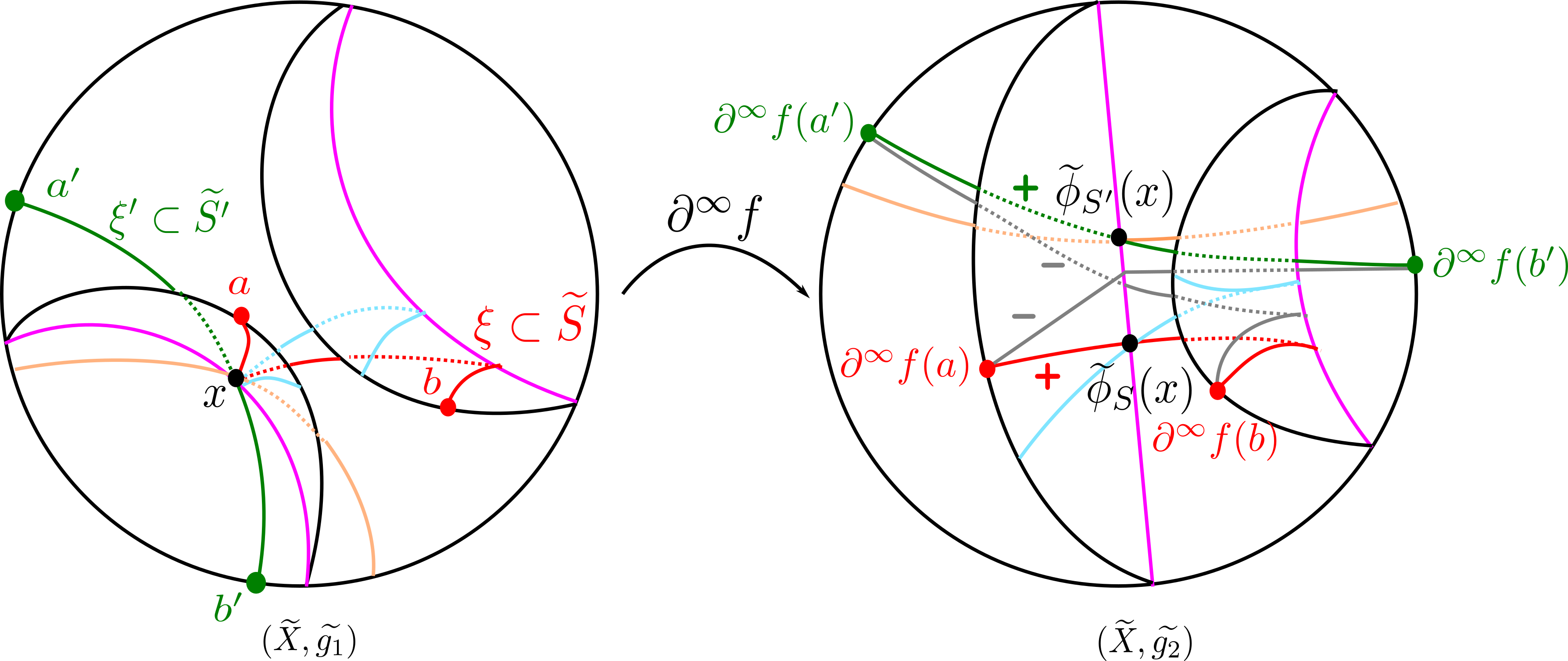}
    \caption{If $\widetilde{\phi}_S(x) \neq \widetilde{\phi}_{S'}(x)$, then $[\partial^{\infty}f(a) \partial^{\infty}f(a) \partial^{\infty}f(b_i) \partial^{\infty}(b'_i)] \neq 0$, but $a$, $a'$, $b$, and $b'$ were chosen so that $[aa'bb'] = 0$, a contradiction since $\partial^{\infty}f$ is Mobius. As before, the pink geodesics are branching.} 
    \label{fig:branchingisom}
\end{figure}

We then conclude that $\partial^{\infty}f$ does not preserve the cross ratio, a contradiction. It then follows that $\widetilde{\phi_S}$ and $\widetilde{\phi_{S'}}$ must pointwise agree on $\widetilde{\gamma}$, as claimed. Finally, $\pi_1$-invariance of Otal's maps gives us our desired result.  
\end{proof}

\begin{lemma} 
\label{chamberagree}
Suppose that $S$ and $S'$ are two closed surfaces in $X$, and suppose $S \cap S' = \{C_i\}_{i = 1}^n$, where each $C_i$ is a chamber. Then for all $C_i \in S \cap S'$, $\phi_S(x) = \phi_{S'}(x)$ where as before, $\phi_S: (S, g_1\vert_{S}) \rightarrow (S, g_2\vert_{S})$ and $\phi_{S'}: (S', g_1\vert_{S'}) \rightarrow (S', g_2\vert_{S'})$ are the isometries constructed by Otal. 
\end{lemma}

\begin{proof} By Lemma \ref{chamberisom}, it suffices to prove $\phi_S(x) = \phi_{S'}(x)$ for $x \in \text{Int}(C_i)$, where $C_i \subset S \cap S'$. Suppose $\widetilde{x}$ is a lift of $x$ inside the interior of $\widetilde{C_i}$, some polygonal lift of $C_i$ in $\widetilde{X}$. Furthermore, suppose $\widetilde{C_i}$ is adjacent to some $\widetilde{\gamma_1}$, a lift of a gluing curves $\gamma \in S \cap S'$. 

We first find $\xi = (w_1, z_1)$ and $\eta = (w_2, z_2)$ on $\widetilde{S}$, some lift of $S$ in $\widetilde{X}$, such that $\xi$ and $\eta$ intersect at $\widetilde{x}$. Recall that  $\widetilde{\phi_S}(\widetilde{x}) = (\partial^{\infty}f(w_1), \partial^{\infty}f(z_1)) \cap (\partial^{\infty}f(w_2), \partial^{\infty}f(z_2))$, where $\partial^{\infty}f$ is the boundary homeomorphism induced by the identity on $X$. Similarly, we choose $\xi' = (w'_1, z'_1)$ and $\eta' = (w'_2, z'_2)$ in a copy of $\widetilde{S'}$ such that $\xi' \cap \eta' = \widetilde{x}$ so $\widetilde{\phi_{S'}}(\widetilde{x}) = (\partial^{\infty}f(w'_1), \partial^{\infty}f(z'_1)) \cap (\partial^{\infty}f(w'_2), \partial^{\infty}f(z'_2))$. We want to use our choice of $\xi$, $\xi'$, $\eta$, and $\eta'$ to show that $\widetilde{{\phi}_{S}}(\widetilde{x}) = \widetilde{{\phi}_{S'}}(\widetilde{x})$. 

In order to construct $\xi$, consider a geodesic arc $\alpha = [\widetilde{x}, \widetilde{\gamma_1}(a)]$ joining $\widetilde{x}$ to some point $\widetilde{\gamma_1}(a)$ on $\widetilde{\gamma_1} \subset \widetilde{X}$. Furthermore, we require that $\widetilde{\gamma_1}(a)$ is chosen so that $\alpha$ lies entirely inside $\widetilde{C_i}$. 
To find $\xi$, extend $\alpha$ in $\widetilde{S}$ via the exponential map to a geodesic lying entirely inside a single lift of $S$, $\widetilde{S}$. Similarly, construct $\xi' \subset \widetilde{S'}$ by extending $\alpha$ to a geodesic lying completely inside $\widetilde{S'}$. Along the same vein, construct $\eta$ by extending $\beta = [\widetilde{x}, \widetilde{\gamma_1}(b)]$ on $\widetilde{S}$, where $\widetilde{\gamma_1}(b)$ is another point on the image of $\widetilde{\gamma_1}$ and $\beta$ also lies entirely inside $\widetilde{C_i}$. Similarly, we can extend $\beta$ in $\widetilde{S'}$ to obtain $\eta'$. See Figure \ref{fig:agree} for an example of a choice of $\alpha$ and $\beta$ given some $\widetilde{x}$. 

We thus have two pairs of geodesics $\xi, \eta \subset \widetilde{S}$, and $\xi', \eta' \subset \widetilde{S'}$. Let $p^{-1}(\{C_i\})$ denote the collection of lifts of chambers in $\{C_i\}_{i = 1}^n = S \cap S'$ in $\widetilde{X}$. Note that $\xi$ and $\xi'$ either agree on a geodesic arc $[\widetilde{\gamma_2}(c), \widetilde{\gamma_1}(a)]$ (where $\widetilde{\gamma_2}$ is a lift of a gluing curve in $S \cap S'$) or a geodesic ray with endpoint $\widetilde{\gamma_1}(a)$ and limit point $y \in \partial^{\infty}\big(p^{-1}(\{C_i\})\big)$. Similarly, $\eta$ and $\eta'$ agree either on a geodesic arc $[\widetilde{\gamma_3}(d), \widetilde{\gamma_1}(b)]$ or a geodesic ray with endpoint $\widetilde{\gamma_1}(b)$ and limit point $z \in \partial^{\infty}\big(p^{-1}(\{C_i\})\big)$. For example, in Figure \ref{fig:agree}, $\xi \subset \widetilde{S}$ and $\xi' \subset \widetilde{S'}$ agree on a geodesic arc $[\widetilde{\gamma}_1(a), \widetilde{\gamma_2}(c)]$ while $\eta \subset \widetilde{S}$ and $\eta' \subset \widetilde{S'}$ agree on a geodesic ray $[\widetilde{\gamma_1}(b), y_2)$. 

\begin{figure}[h!]
    \centering    \includegraphics[width=\textwidth]{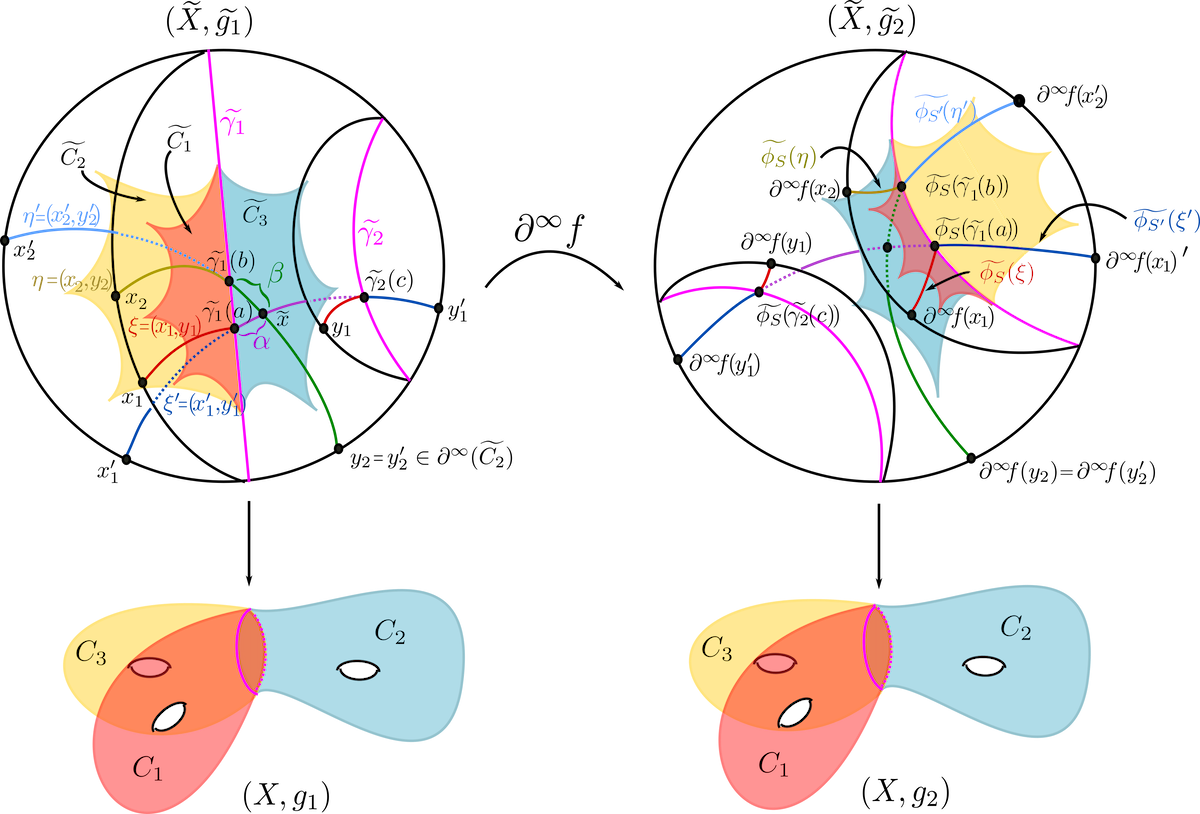}
    \caption{We extend $\alpha = (\widetilde{x}, \widetilde{\gamma_1}(a))$ in the picture on the left to $\xi \subset \widetilde{S}$ and $\xi' \subset \widetilde{S'}$. Similarly, we extend $\beta = (\widetilde{x}, \widetilde{\gamma_1}(b))$ where $b \neq a$ to $\eta \subset \widetilde{S}$ and $\eta' \subset \widetilde{S'}$. While $\xi$ and $\xi'$ agree on a geodesic arc $[\widetilde{\gamma_2}(c), \widetilde{\gamma_1}(a)]$, $\eta$ and $\eta'$ agree on a geodesic ray $[\widetilde{\gamma_1}(b), y)$.}
    \label{fig:agree}
\end{figure}

Suppose that $\xi$ and $\xi'$ agree on a geodesic arc $[\widetilde{\gamma_2}(c), \widetilde{\gamma_1}(a)] \subset p^{-1}(\{C_i\})$. By Proposition \ref{chamberisom}, we know that we can define points $x_1 := \widetilde{\phi_S}(\widetilde{\gamma_1}(a)) = \widetilde{\phi_{S'}}(\widetilde{\gamma_1}(a))$ and $x_2 := \widetilde{\phi_S}(\widetilde{\gamma_2}(c)) = \widetilde{\phi_{S'}}(\widetilde{\gamma_2}(c))$. Then it follows that $\widetilde{\phi_S}(\xi)$ and $\widetilde{\phi_{S'}}(\xi')$ agree on the unique geodesic arc between $x_1$ and $x_2$. In other words, 
$$\widetilde{\phi_S}(\xi \cap \xi') = \widetilde{\phi_S}([\widetilde{\gamma_1}(a), \widetilde{\gamma_2}(c)]) = [x_1, x_2] = \widetilde{\phi_{S'}}([\widetilde{\gamma_1}(a), \widetilde{\gamma_2}(c)]) = \widetilde{\phi_{S'}}(\xi \cap \xi').$$

Suppose on the other hand that $\xi$ and $\xi'$ agree on a geodesic ray $[\widetilde{\gamma_1}(a), y)$. Again, we can set $x_1 := \widetilde{\phi_S}(\widetilde{\gamma_1}(a)) = \widetilde{\phi_{S'}}(\widetilde{\gamma_1}(a))$. Note that we can define $y_1 := \partial^{\infty}\widetilde{\phi_S}(y) = \partial^{\infty}\widetilde{\phi_{S'}}(y)$ since in the proof of Lemma \ref{chamberisom}, we determined that $\partial^{\infty}\widetilde{\phi_S}\vert_{\widetilde{S \cap S'}} = \partial^{\infty}f\vert_{\widetilde{S \cap S'}} = \partial^{\infty}\widetilde{\phi_{S'}}\vert_{\widetilde{S \cap S'}}$. Since the geodesic between $x_1$ and $y_1$ is unique in $\widetilde{X} \cup \partial^{\infty}(\widetilde{X})$, it then follows that $\widetilde{\phi_S}(\xi)$ and $\widetilde{\phi_{S'}}(\xi')$ necessarily agree on a geodesic ray with endpoints $x_1$ and $y_1$. Thus,
$$\widetilde{\phi_S}(\xi \cap \xi') = \widetilde{\phi_S}([\widetilde{\gamma_1}(a), y)) = [x_1, y_1) = \widetilde{\phi_{S'}}([\widetilde{\gamma_1}(a), y)) = \widetilde{\phi_{S'}}(\xi \cap \xi').$$

It then follows that $\widetilde{\phi_{S}}(\xi \cap \xi')$ and $\widetilde{\phi_{S'}}(\xi \cap \xi')$ agree, and the same argument holds for $\widetilde{\phi_{S}}(\eta \cap \eta')$ and $\widetilde{\phi_{S'}}(\eta \cap \eta')$. We now show that $\widetilde{\phi}(\xi \cap \xi') := \widetilde{\phi_S}(\xi \cap \xi') = \widetilde{\phi_{S'}}(\xi \cap \xi')$ and $\widetilde{\phi}(\eta \cap \eta') := \widetilde{\phi_S}(\eta \cap \eta') = \widetilde{\phi_{S'}}(\eta \cap \eta')$ intersect at a point $\widetilde{x} \in 
\widetilde{\phi_S}(\widetilde{S} \cap \widetilde{S'}) = \widetilde{\phi_{S'}}(\widetilde{S} \cap \widetilde{S'})$ by doing some case work. 

\begin{figure}[H]
    \centering    \includegraphics[width=0.95\textwidth]{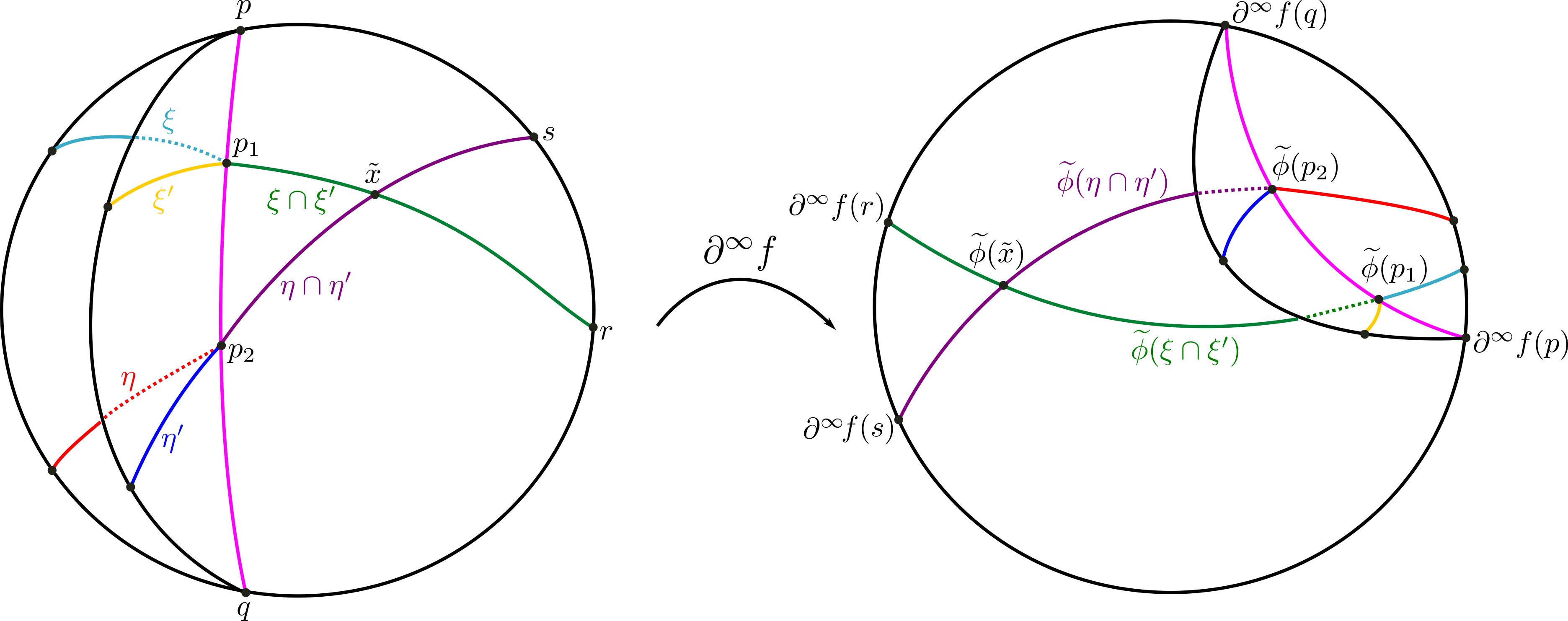}
    \caption{An illustration of Case 1; $\xi \cap \xi'$ and $\eta \cap \eta'$ are both geodesic rays.}
    \label{fig:intersect1}
\end{figure}

\begin{figure}[H]
    \centering
    \includegraphics[width=0.95\textwidth]{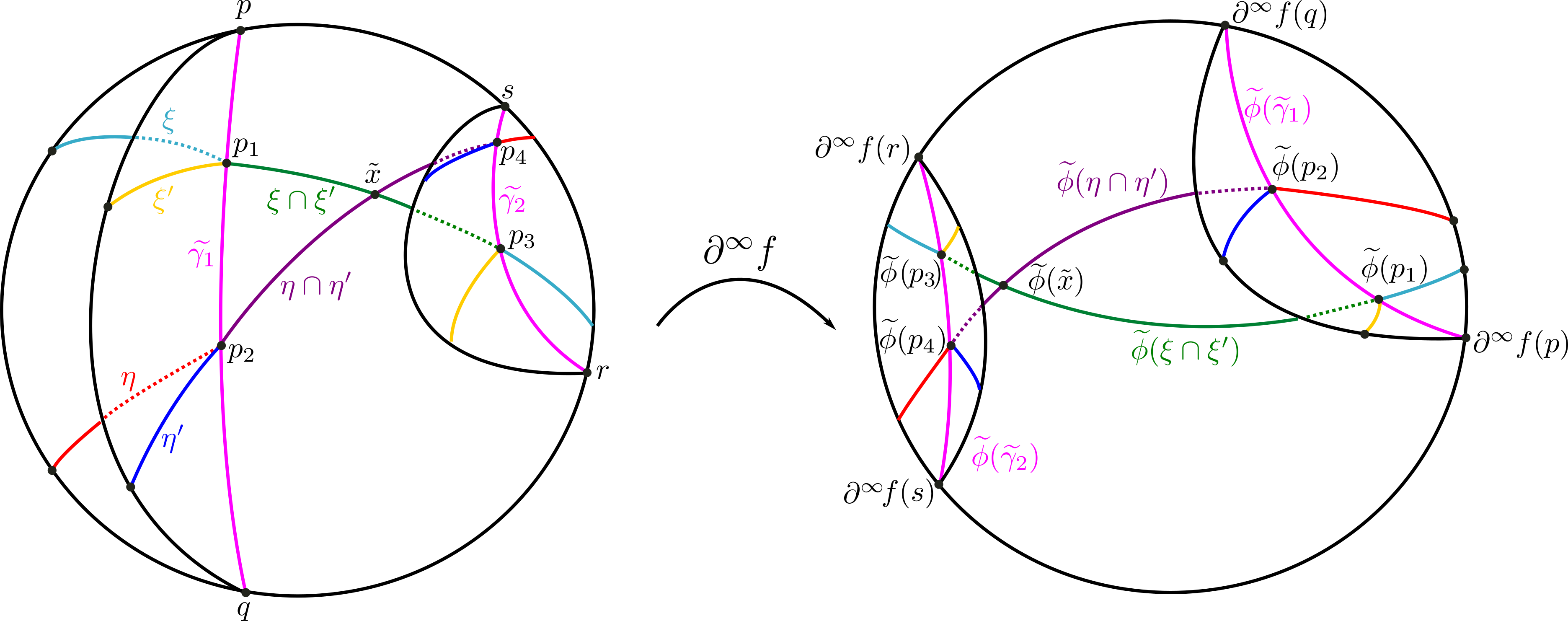}
    \caption{An example of Case 2; $\xi \cap \xi'$ and $\eta \cap \eta'$ are both geodesic arcs.}
    \label{fig:intersect2}
\end{figure}

\begin{figure}[H]
    \centering
    \includegraphics[width=0.95\textwidth]{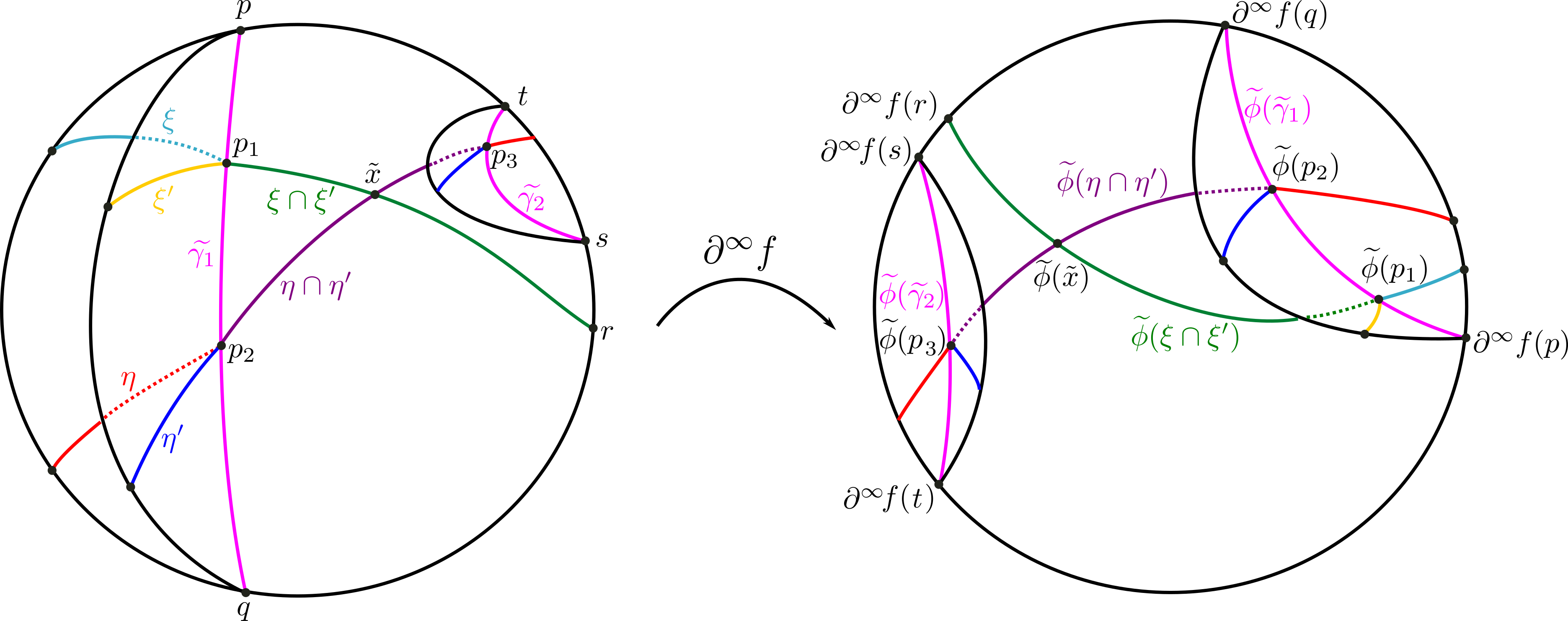}
    \caption{An example of Case 3; $\xi \cap \eta$ is a geodesic ray and $\xi' \cap \eta'$ is a geodesic arc.}
    \label{fig:intersect3}
\end{figure}

\begin{enumerate}
    \item \textbf{Case One: $\xi \cap \xi'$ and $\eta \cap \eta'$ are both geodesic rays.} Suppose $\xi \cap \xi' = [p_1, r)$ and $\eta \cap \eta' = [p_2, s)$, where $p_1 = \widetilde{\gamma_1}(a)$ and $p_2 = \widetilde{\gamma_1}(b)$ lie on a branching geodesic $\widetilde{\gamma_1}$ and $r, s \in \partial^{\infty}\big(p^{-1}(\{C_i\})\big)$. Additionally, suppose $\widetilde{\gamma_1}$ has endpoints $p$ and $q$, $p = \lim\limits_{t \rightarrow -\infty}\widetilde{\gamma_1}(t)$, and $q = \lim\limits_{t \rightarrow \infty}\widetilde{\gamma_1}(t)$. Consider an apartment $A$ in $\widetilde{X}$ that contains $\widetilde{\gamma_1}$, $r$, and $s$ (e.g. $\widetilde{S}$ or $\widetilde{S'}$). By construction, if without loss of generality $a < b$ (so that $p_1$ is ``closer" to $p$ and $p_2$ is ``closer" to $q$), the order of points, according to some fixed orientation on $\partial^{\infty}(A)$, will be $p, q, r, s$, and $\partial^{\infty}f\vert_A$ will preserve that order. Furthermore, since the orientation of $\widetilde{\gamma_1}$ will be preserved under $\partial^{\infty}f\vert_A$, it follows that the geodesic rays $[\widetilde{\phi}(p_2), \partial^{\infty}f(s)) = [\widetilde{\phi}(p_2), \partial^{\infty}\widetilde{\phi}(s))$ and $[\widetilde{\phi}(p_1), \partial^{\infty}f(r)) = [\widetilde{\phi}(p_1), \partial^{\infty}\widetilde{\phi}(r))$ will necessarily intersect (see Figure \ref{fig:intersect1}). 
    
    \item \textbf{Case Two: $\xi \cap \xi'$ and $\eta \cap \eta'$ are both geodesic arcs.} Suppose $\xi \cap \xi' = [p_1, p_3]$ and $\eta \cap \eta' = [p_2, p_4]$. Let $p$ and $q$ be defined as before in case one. As before, suppose $\widetilde{\gamma_1}(a) = p_1 \in \xi \cap \xi'$ and $\widetilde{\gamma_1}(b) = p_2 \in \eta \cap \eta'$ where $a < b$. Additionally, suppose that $\widetilde{\gamma_2}(c) = p_3$ and $\widetilde{\gamma_2}(d) = p_4$ where $d < c$, and $\lim\limits_{t \rightarrow -\infty}\widetilde{\gamma_2}(t) = s$ and $\lim\limits_{t \rightarrow \infty} \widetilde{\gamma_2}(t) = r$. Again, one can find some apartment $A $ (e.g. $\widetilde{S}$ or $\widetilde{S'}$) that contains $\gamma_1$ and $\gamma_2$, and as before, we have that $\partial^{\infty}f\vert_A$ is a circle homeomorphism which preserves the cyclic order of the four-tuple $[p, q, r, s]$. Since the orientation of $\widetilde{\gamma_1}$ and $\widetilde{\gamma_2}$ is, as before, preserved, it follows that $[\widetilde{\phi}(p_1), \widetilde{\phi}(p_3)]$ and $[\widetilde{\phi}(p_2), \widetilde{\phi}(p_4)]$ intersect (see Figure \ref{fig:intersect2}).  
    
    \item \textbf{Case Three: $\xi \cap \xi'$ or $\eta \cap \eta'$ is a geodesic ray, and the other is a geodesic arc.} For the last case, suppose without loss of generality that $\xi \cap \xi'$ is a geodesic ray $[p_1, r)$ and $\eta \cap \eta'$ is a geodesic arc $[p_2, p_3]$. Again, suppose $p_1 = \widetilde{\gamma_1}(a)$ and $p_2 = \widetilde{\gamma_1}(b)$ for some $a < b$ where $\widetilde{\gamma_1}$ has endpoints $p = \lim\limits_{t \rightarrow -\infty}\widetilde{\gamma_1}(t)$ and $q = \lim\limits_{t \rightarrow \infty} \widetilde{\gamma_1}(t)$. Suppose $p_3 = \widetilde{\gamma_2}(c)$ where $c \in \mathbb{R}$ and the endpoints of $\widetilde{\gamma_2}$ are $s$ and $t$. Again, choose some apartment $A$ (e.g. $\widetilde{S}$ or $\widetilde{S'}$) that contains all five points ($p, q, r, s,$ and $t$). Then, after fixing an orientation, the points will be ordered $p, q, r, s, t$, and their order will be preserved under $\partial^{\infty}f\vert_A$. Again, since the orientation of $\widetilde{\gamma_1}$ is fixed, it follows that $[\widetilde{\phi}(p_1), \partial^{\infty}f(r)) = [\widetilde{\phi}(p_1), \partial^{\infty}\widetilde{\phi}(r))$ and $[\widetilde{\phi}(p_2), \widetilde{\phi}(p_3)]$ will intersect. (See Figure \ref{fig:intersect3}.)
    
\end{enumerate}

In conclusion, $\widetilde{\phi_S}(\xi \cap \xi') = \widetilde{\phi_{S'}}(\xi \cap \xi')$ and $\widetilde{\phi_{S}}(\eta \cap \eta') = \widetilde{\phi_{S'}}(\eta \cap \eta')$ intersect at some point in $\widetilde{\phi_S}\big(p^{-1}(\{C_i\})\big) = \widetilde{\phi_{S'}}\big(\{p^{-1}(\{C_i\}\big)$. Since $\widetilde{\phi_S}(\xi \cap \eta) = \widetilde{\phi_S}(\widetilde{x})$, $\widetilde{\phi_{S'}}(\xi' \cap \eta') = \widetilde{\phi_{S'}}(\widetilde{x})$ and $\widetilde{\phi_S}$ and $\widetilde{\phi_{S'}}$ are injective, it follows that:
$$\widetilde{\phi_{S'}}(\widetilde{x}) = \widetilde{\phi_{S'}}(\xi \cap \xi' \cap \eta \cap \eta') = \widetilde{\phi_{S'}}(\xi \cap \xi') \cap \widetilde{\phi_{S'}}(\eta \cap \eta') = \widetilde{\phi_{S}}(\xi \cap \xi') \cap \widetilde{\phi_{S}}(\eta \cap \eta') = \widetilde{\phi_S}(\xi \cap \xi' \cap \eta \cap \eta') = \widetilde{\phi}_{S}(\widetilde{x}).$$
By the $\pi_1(S)$ and $\pi_1(S')$-equivariance of $\phi_S$ and $\phi_{S'}$ respectively, it follows that $\phi_S(x) = \phi_{S'}(x)$. 
\end{proof}

The proof of Proposition \ref{main1} then follows easily: \\

\noindent 
\textit{Proof of Proposition \ref{main1}}. As before, consider a minimal collection of closed surfaces $\mathcal{S}$ that covers $X$. Consider any arbitrary $S, S' \in \mathcal{S}$. By Lemmas \ref{chamberisom} and \ref{chamberagree}, $\phi_S$ and $\phi_{S'}$ pointwise agree on $S \cap S'$. Since $X$ is complete, convex, and compact, by Lemma \ref{isom}, one can thus patch the collection of isometries $\{\phi_S\}_{S \in \mathcal{S}}$ together to obtain a global isometry $\phi: (X, g_1) \rightarrow (X, g_2)$. 
\qed 
\section{The General Case}\label{generalcase}

We now tackle the general case of \Cref{main}. We will begin by showing that simple (\textit{not} necessarily thick), NPC surface amalgams are NPC cube complexes. Then, using some powerful machinery mentioned in \Cref{sec:qcerf}, we deduce that simple, negatively curved surface amalgams are QCERF. Next, we show that simple, thick, negatively curved surface amalgams can be covered by finitely many immersed closed surfaces. Finally, since surface subgroups of $\pi_1(X)$ are quasiconvex, we use Scott's theorem (\Cref{thm:scott}) to promote the immersed closed surfaces to embedded closed surfaces in finite-sheeted covers, allowing us to reduce to the base case.

As promised, we first prove the following:

\begin{lemma}\label{lemma:NPC}
Let $X$ be a simple, NPC surface amalgam. Then $X$ is a NPC cube complex. 
\end{lemma}

\begin{proof} We will realize $X$ as a cube complex and then apply the combinatorial link condition to show $X$ is NPC. First, we assume $X$ is negatively curved. Since the gluing geodesics of $X$ are simple and disjoint, there exists a pants decomposition $\mathcal{P}$ of $X$ such that the gluing curves are contained in the set of pants curves in $\mathcal{P}$. We first realize each pair of pants $P \in \mathcal{P}$ as cube complex: 

\begin{figure}[H]
    \centering
\includegraphics[width=0.5\textwidth]{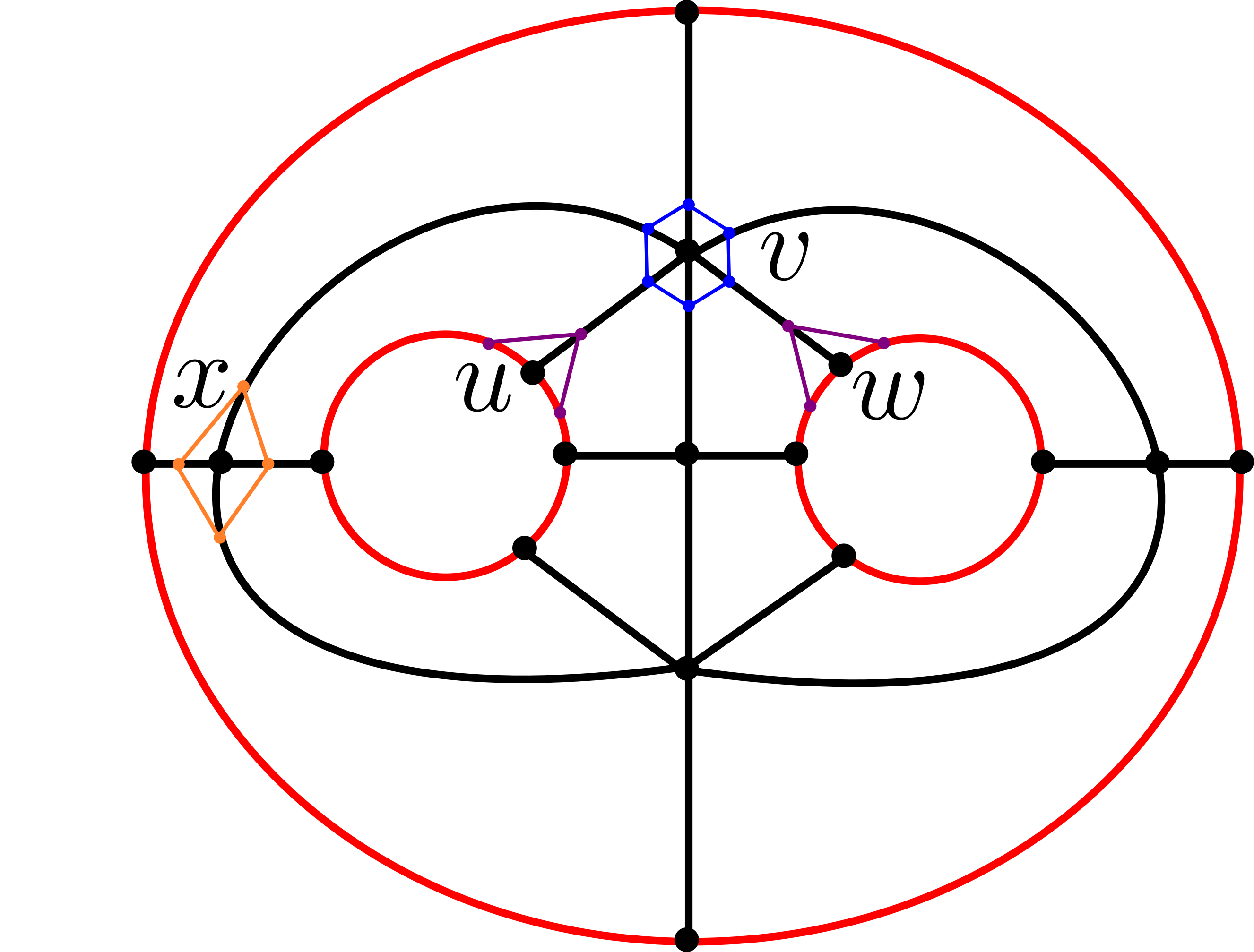}
    \caption{A pair of pants as a cube complex. The boundary components are red. The link $lk(v)$ is shown in blue, $lk(u)$ and $lk(w)$ are shown in purple, and $lk(x)$ is in orange.}
    \label{fig:pants}
\end{figure}

Note that each boundary component of $P$ contains four vertices. Given a collection of boundary components $\{b_i\}$ identified in $X$, one can subdivide each $b_i$ so that the vertices are aligned when the $b_i$ are glued together. The result is a cube complex structure for $X$.  

We claim that for every vertex $v$ in the cube complex structure of $X$, $lk(v)$ is triangle-free. Indeed, we list all the possibilities: 
\begin{enumerate} 
\item If $v$ is one of two vertices in each pair of pants adjacent to six edges, $lk(v)$ will be a hexagon (see blue link in \Cref{fig:pants});
\item If $v$ is one of three vertices in each pair of pants adjacent to four edges, then $lk(v)$ will be a square (see orange link in \Cref{fig:pants});
\item If $v$ lies on a pants curve not attached to any other pants curve, then its link is a graph with three vertices and two edges (see purple links in \Cref{fig:pants});
\item Finally, if $v$ lies on a pants curve glued to at least one other pants curve, then $lk(v)$ is a union of squares all sharing two vertices lying on boundary components of the pants decompositions. For instance, if $u$ and $w$ are identified in \Cref{fig:pants}, then the union of the purple links would form a square.
\end{enumerate} 

In all four cases, $lk(v)$ is triangle-free, and thus would not span a simplex of dimension greater than $1$. Then all the links are flag complexes, as desired.

The generalization to the NPC setting follows easily; one can naturally represent each cylindrical chamber with a square with one pair of opposite sides identified. Subdivide the unidentified sides so that each side has four vertices, including the pair of identified vertices on the corners of the original square. As before, attach boundary components together so that vertices are identified with other vertices, giving the surface amalgam a cube complex structure. Now, the links of each vertex on the cylinders can only fall under case (3) or (4) from the aforementioned list, thus satisfying the combinatorial link condition. 
\end{proof} 

Next, we use \Cref{lemma:NPC} to show that the fundamental group of a simple, negatively curved surface amalgam is QCERF: 

\begin{lemma}\label{lemma:qvh2} Let $X$ be a simple, negatively curved surface amalgam. Then $\pi_1(X)$ is QCERF. 
\end{lemma}

\begin{proof} We collect a few useful facts and observations:

\begin{enumerate}
\item \textit{$\mathbb{F}_n \in \mathcal{QVH}$ for all $n \in \mathbb{N}$.} Recall that by \Cref{thm:qvh}, a torsion-free Gromov hyperbolic group is in $\mathcal{QVH}$ if and only if it is virtually special. From this, we can naturally deduce that any finite rank free group, which is subgroup separable by a classic theorem by Hall (\cite{hall}), is QCERF and therefore virtually special, and in $\mathcal{QVH}$.
\item \textit{Infinite cyclic subgroups of $\pi_1(X)$ are quasiconvex}. This follows from the fact that if a group $G$ (such as $\pi_1(X)$) is Gromov hyperbolic, every infinite cyclic subgroup of $G$ is quasiconvex. 
\item \textit{The fundamental group of a surface amalgam is torsion-free.} Due to Proposition 2 of \cite{serre}, every finite order element of an amalgamated product $A \ast_C B$ (resp. HNN extension $A \ast_C$) is conjugate to an element of $A$ or $B$ (resp. $A$). Since free groups are torsion free, one can inductively show that the fundamental group of any surface amalgam is torsion-free. 
\end{enumerate} 

By definition of $\mathcal{QVH}$, (1) and (2), $\pi_1(X) \in \mathcal{QVH}$. By (3), Gromov hyperbolicity of $\pi_1(X)$, and \Cref{thm:qvh}, $\pi_1(X)$ virtually special. Thus, by \Cref{lemma:NPC} and \Cref{thm:hw}, it follows that $\pi_1(X)$ is QCERF as well. 
\end{proof} 

We briefly mention the following corollary, which could be of independent interest: 

\begin{corollary}
    Let $(X, g)$ be a simple, negatively curved surface amalgam. Then $\pi_1(X)$ is residually finite.
\end{corollary}

Finally, we show that every chamber is embedded in an immersed closed surface in $X$: 

\begin{lemma}\label{thm:immersion} Let $(X, g)$ be a simple, thick, negatively curved surface amalgam. Then each chamber in $(X, g)$ can be included in an embedded closed surface in some finite-sheeted cover $(\widehat{X}, \widehat{g})$ of $(X, g)$. 
\end{lemma}

\begin{proof} We claim that each surface amalgam can be covered by a disjoint union of immersed closed surfaces. We first construct this collection of closed surfaces. 

Take two copies of each chamber in $C_i \subset X$ to obtain two identical collections of chambers $\{C_i^1\} \cup \{C_i^2\}$. Consider a gluing geodesic $\gamma \subset X$, and suppose that $\{b_k\} \subset X$ is the collection of boundary components identified with $\gamma$ in $X$ indexed by some collection of natural numbers $k = 1, 2, ..., N$. We will label the two copies of each $b_k$ with $b_k^1$ and $b_k^2$. For $1 \leq k \leq N$, identify $b_k^1$ with $b_{k + 1}^2$, where $k + 1$ is taken modulo $N$. Repeat the same process for every gluing curve in $X$ so that each boundary component of each chamber is glued to exactly one other boundary component. Due to the thickness assumption, the result is a (possibly disconnected) collection of finitely many closed surfaces $\{N_s\}$. An example of the process is illustrated in \Cref{fig:immersion}. 

\begin{figure}[H]
    \centering
\includegraphics[width=\textwidth]{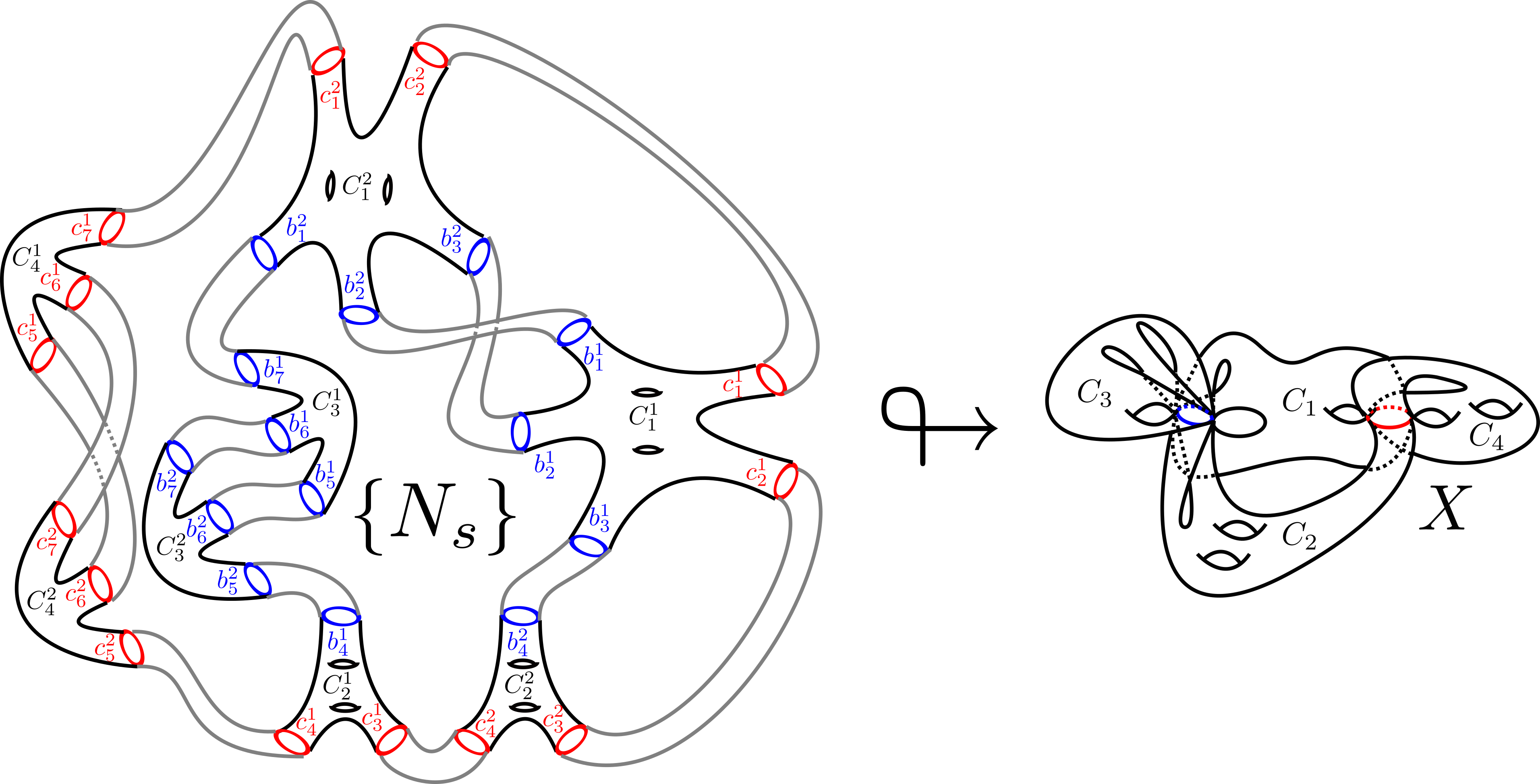}
    \caption{An immersion of a closed surface into a surface amalgam $X$ using the construction from the proof of \Cref{thm:immersion}.}
    \label{fig:immersion}
\end{figure}

We now show that there is an immersion $f: \{N_s\} \looparrowright X$, given by the projection of $\{N_s\}$ back onto $X$. More specifically, if $y \in \{N_s\} \setminus (\{b_k^1\} \cup \{b_k^2\})$, then it projects to the interior of a chamber in $X$. Otherwise, $y$ projects to a point on a gluing geodesic in $X$. 

Usually, an immersion between manifolds is defined as a map with a locally injective derivative. An analogous definition holds for surface amalgams, but we need to define the local unit tangent bundle carefully. Let $S'_xC$ be the open hemisphere of unit tangent vectors based at $x$ pointing inside a chamber $C$ containing or adjacent to $x$. Then define $S'_x := \bigcup\limits_{i = 1}^{m} S_x'C_i$, the collection of the $m$ open hemispheres of unit tangent vectors based at $x$ and pointing into the $m$ chambers adjacent to or containing $x$ (see Figure \ref{fig:vector}). Note that if $x \in X$ lies on the interior of a chamber, $S'_x(X)$ resembles the usual circle of unit tangent vectors based at a point on a surface.

\begin{figure}[h!]
    \centering
    \includegraphics[width=0.6\textwidth]{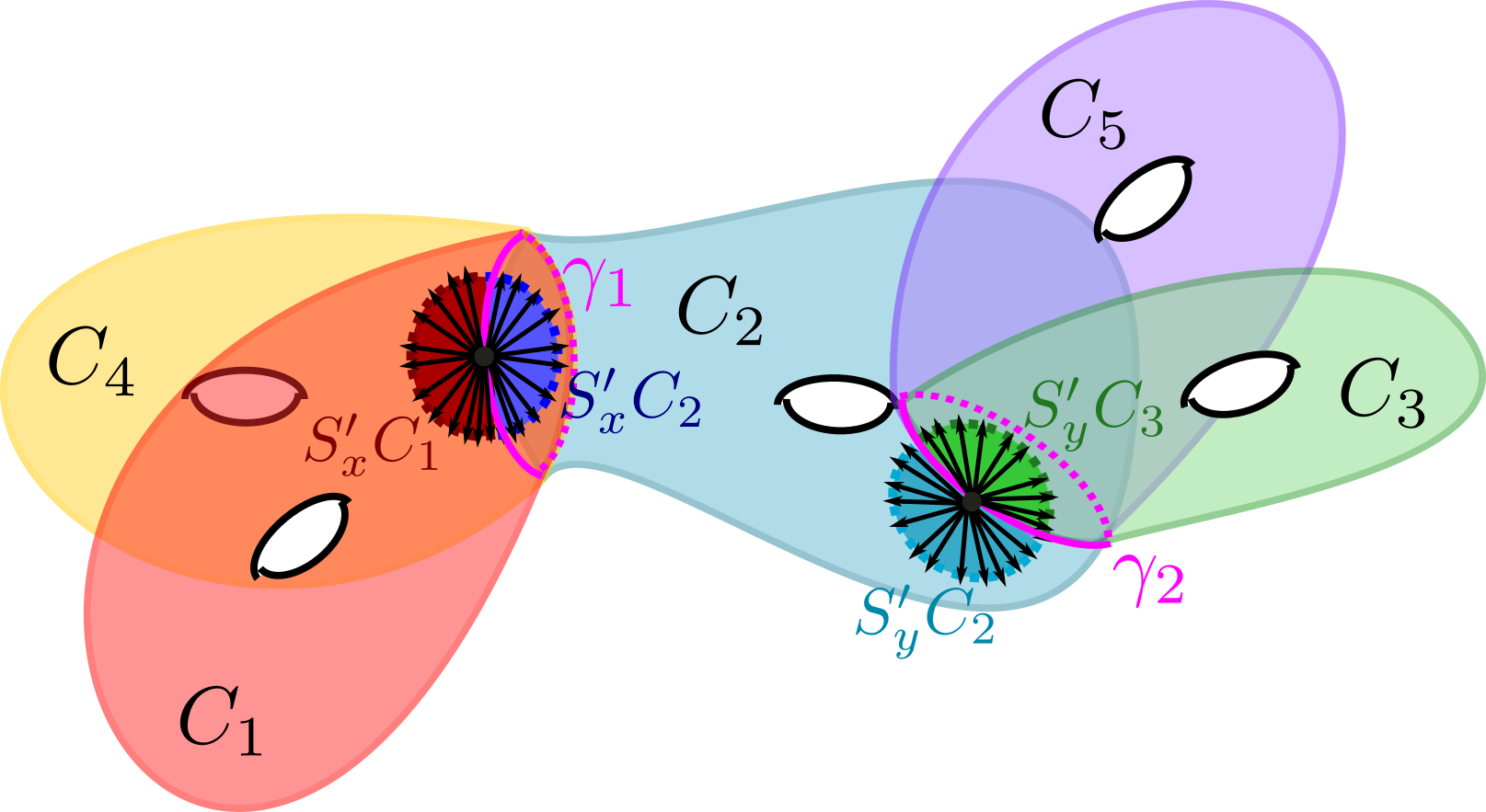}
    \caption{Given two points $x \in \gamma_1$ and $y \in \gamma_2$, we illustrate examples of hemispheres of unit tangent vectors based at $x$ and $y$. Here, $S'_x = S'_xC_1 \cup S'_xC_2 \cup S'_xC_4$ and $S'_y = S'_yC_2 \cup S'_yC_3 \cup S'_yC_5$.}
    \label{fig:vector}
\end{figure}

With this definition in mind, note that for all $y \in \{N_s\}$, $Df_y: T_y (\{N_s\}) \rightarrow S'_{f(y)}(X)$ is injective. We point out that injectivity comes from the fact that there do not exist $b_k^1$ and $b_k^2$ which are glued together in $\{ N_s\}$. In conclusion, each chamber $C$ sits inside an immersed closed surface in $X$, which we will call $N_C$. 

The preimage of $N_C$ in the universal cover $\widetilde{X}$ under the covering map is a convex set, as it is a union of polygons with geodesic boundary. As a consequence, $\pi_1(N_C)$ is a quasiconvex subgroup of $\pi_1(X)$. Since $\pi_1(X)$ is QCERF, there is some finite-sheeted cover $\widehat{X}$ of $X$ such that $N_C \subset \widehat{X_C}$ embeds as a closed surface due to \Cref{thm:scott}. Given a collection of finite-sheeted covers $\{\widehat{X_C}\}$, there is a common finite-sheeted cover $\widehat{X}$ such that every chamber $C$ embeds in a closed surface, as the intersection of a finite-index subgroup is still a finite-index subgroup. This proves the lemma.  

\end{proof} 

\subsection{Proof of Theorem \ref{main}}

We now have all the ingredients needed for proving Theorem \ref{main}.  

\begin{proof} Suppose $(X, g_1)$ and $(X, g_2)$ have the same marked length spectra. Suppose every chamber in $X$ can be included into a closed surface. Then by Proposition \ref{main1}, $(X, g_1)$ and $(X, g_2)$ are isometric via a map isotopic to identity, $\phi: (X, g_1) \rightarrow (X, g_2)$. 

Otherwise, by \Cref{thm:immersion}, there exists some finite-sheeted cover $\widehat{X}$ of $X$ such that every chamber is embedded in $\widehat{X}$. Then by Proposition \ref{main1}, there exists an isometry isotopic to the identity $\widehat{\phi}: (\widehat{X}, \widehat{g_1}) \rightarrow (\widehat{X}, \widehat{g_2})$. Recall that $\widehat{\phi}$ is constructed by projecting the $\pi_1(\widehat{X})$-equivariant isometry $\widetilde{\phi}: (\widetilde{X}, \widetilde{g_1}) \rightarrow (\widetilde{X}, \widetilde{g_2})$ between copies of the universal cover $\widetilde{X}$ of both $X$ and $\widehat{X}$. By construction, $\widetilde{\phi}$ is also $\pi_1(X)$-equivariant, so there exists some isometry $\phi: (X, g_1) \rightarrow (X, g_2)$ isotopic to the identity as well, as desired.
\end{proof}

\section{Appendix} 

\label{ergodic}

The main purpose of this appendix is to prove the following statement, which is essential to the proof of Lemma \ref{approx}: 

\begin{proposition}
\label{ergodicity} 
The geodesic flow map on (Gromov) hyperbolic P-manifolds is ergodic with respect to the Bowen-Margulis measure.
\end{proposition}

Notice that for the purposes of the proof of Lemma \ref{approx}, it does not really matter what measure the geodesic flow map is ergodic with respect to. We remark that in the case of symmetric spaces, the Bowen-Margulis measure coincides with the Liouville measure.

The ideas and terminology of this section are mostly taken from \cite{kai}, which we give a summary of for the convenience of the reader. We remark that we can also deduce ergodicity of the geodesic flow map from general theory in the setting of CAT(0) spaces with rank-one axes developed in \cite{ricks}, but for simplicity, we stick to the ideas in \cite{kai}. 

For the rest of this subsection, let $(X, d)$ be a proper, connected Gromov hyperbolic space under a proper (preimages of compact sets are compact), nonelementary, isometric action by a group $\Gamma$ and $\partial^{\infty}(X, d)$ be its visual boundary, the set of equivalence classes of asymptotic geodesic rays (see Definition \ref{boundary} for details). Then the \textit{(Gromov) hyperbolic compactification} $\overline{X}$ of $X$ may be defined as $\overline{X} = X \cup \partial^{\infty}(X)$. 

Following \cite{kai}, we impose two additional assumptions:

\begin{assumption}[Uniqueness of geodesics]
    \label{U} For all $x_1, x_2 \in \overline{X}$, there exists a unique geodesic $[x_1, x_2]$ joining them. 
\end{assumption}

\begin{assumption}[Existence of convergent geodesic rays]
    \label{C} If $\alpha_1$ and $\alpha_2$ are asymptotic (e.g. they lie within bounded distance of each other) geodesic rays, then there exists some $c > 0$ such that: $$\lim\limits_{t \rightarrow \infty} d\big(\alpha_1(t), \alpha_2(t + c)\big) = 0.$$ 
\end{assumption}

We remark that both assumptions hold if $(X, d)$ is locally CAT(-1). 

\subsection{Patterson-Sullivan Measures} 

We now define a family of measures $\{\mu_p\}_{p \in X}$ on $\partial^{\infty}(X)$ which are defined for every point $p \in X$. Intuitively, such measures, called \textit{conformal densities}, measure the proportion of elements of $\Gamma x = \{\gamma x \vert \gamma \in \Gamma\}$ that land within a specified subset of $\partial^{\infty}(X)$. 

Given a point $a \in \partial^{\infty}(X)$ and $x, y \in X$, we set: $$B_a(x, y) = \lim\limits_{t \rightarrow \infty} d(x, r(t)) - d(y, r(t))$$
where $r: \mathbb{R} \rightarrow X$ is a geodesic ray with endpoint $a$. Sometimes in the literature, $B_a(x, y)$ is called the \textit{horospherical distance between $x$ and $y$}, and $B_a$, which does not depend on the choice of $r(t)$, is sometimes known as the Busmann cocycle function.  

\begin{definition}\label{conformal}
    A family $\{\mu_p\}_{p \in X}$ of finite Borel (Radon) measures on $\partial^{\infty}(X)$ is called a \textit{conformal density of dimension $\delta$} if: 
    \begin{enumerate}
    \item For all $\gamma \in \Gamma$, $\gamma_{\ast} \mu_p = \mu_{\gamma \cdot p}$ ($\mu_p$ is $\Gamma$-invariant);
    \item For all $p, q \in X$ and $a \in \partial^{\infty}(X)$, $\mu_p$ and $\mu_q$ are equivalent with Radon-Nikodym derivative $$\frac{d\mu_q}{d\mu_p}(a) = e^{-\delta B_a(x, y)}.$$
    \end{enumerate}
\end{definition}

One can define a \textit{quasiconformal density of dimension $\delta$} by relaxing the second condition in Definition \ref{conformal} to say for all $p, q \in X$ and $a \in \partial^{\infty}(X)$, there exists some $C \geq 1$ such that:
$$\frac{1}{C}e^{-\delta B_a(x, y)} \leq \dfrac{d\mu_p}{d\mu_q}(a) \leq Ce^{-\delta B_a(x, y)}.$$

The \textit{Poincare series} associated to $\Gamma$ is the series $P(x, s) = \sum\limits_{\gamma \in \Gamma} e^{sd(x, \gamma.x)}$, where $x \in X$ and $s \in \mathbb{R}$. There is some $\delta_{\Gamma} \in \mathbb{R}$ aptly named the \textit{critical exponent of $\Gamma$}. The series converges if $s < \delta_{\Gamma}$ and diverges when $s > \delta_{\Gamma}$; if $s = \delta_{\Gamma}$, the series could either converge or diverge (see Proposition 5.3 of \cite{coo}). Note that if $X$ is a proper geodesic space and $X/\Gamma$ is compact, then $\delta_{\Gamma}$ is finite (see Proposition 1.7 of \cite{burger}). Furthermore, if the action of $\Gamma$ on $X$ is non-elementary, then the critical exponent is nonzero. Thus, in the case where $X$ is a simple, thick P-manifold and $\Gamma = \pi_1(X)$, $\delta_{\Gamma}$ is nonzero and finite. 

While the existence of conformal densities is not guaranteed for a Gromov hyperbolic space $(X, d)$, one can guarantee the existence of quasiconformal densities given that the critical exponent $\delta_{\Gamma}$ is a finite positive number (as is the case when $X$ is a simple, thick P-manifold), summarized in the following theorem. The \textit{limit set of $\Gamma$}, the set of accumulation points in $\partial^{\infty}(X)$ of $\Gamma x$ for some (any) $x \in X$, is denoted by $\Lambda_{\Gamma}$.

\begin{theorem}[\cite{coo}, Theorem 5.4]\label{qc} 
Suppose a group $\Gamma$ acts properly discontinuously via isometries on $(X, d)$, a proper (Gromov) hyperbolic metric space (and therefore $\delta_{\Gamma} \in (0, \infty)$). Then there exists a quasiconformal density of dimension $\delta_{\Gamma}$ supported on $\Lambda_{\Gamma}$. 
\end{theorem}

For every point $x \in X$, Cooernart constructs a \textit{Patterson-Sullivan measure} $\mu_x$ from the quasiconformal density of dimension $\delta_{\Gamma}$. Its construction depends on whether $P(x, s)$ converges or diverges at $s = \delta_{\Gamma}$; we refer the reader to the proof of Theoerem 5.4 in \cite{coo} for details. In particular, if $\Gamma$ acts convex cocompactly on $X$, then $P(x, s)$ diverges when $s = \delta_{\Gamma}$ (see Corollary 7.3 of \cite{coo}) and the Patterson-Sullivan measure is defined as the weak limit of probability measures:

$$\mu_x = \lim\limits_{n \rightarrow \infty} \frac{1}{\sum\limits_{\gamma \in B_n}e^{-\delta_{\Gamma}d(x, \gamma.x)}}\sum\limits_{\gamma \in B_n} e^{-\delta_{\Gamma}d(x, \gamma.x)}\text{Dirac}_{\gamma.x}$$
where $B_n = \{\gamma \in \Gamma \vert \gamma.o \in B(o, n)\}$ for some fixed $o \in X$. 
Furthermore, he proves that $\text{supp}(\mu_x) = \Lambda_{\Gamma} \subseteq \partial^{\infty}(X)$.

Next, one can use the Patterson-Sullivan measures to define measures on $\partial^{\infty}(X) \times \partial^{\infty}(X) \setminus \Delta$.

\subsection{Bowen-Margulis Measures}

Given $\mu_x$ with density $\delta$, Kaimanovich defines a $\Gamma$-invariant measure on $\partial^{\infty}(X) \times \partial^{\infty}(X) \setminus \Delta$ (see Section 2.4.1 in \cite{kai}): 
$$d\nu_x(a, b) = d\mu_x(a)d\mu_x(b)(e^{{\langle a, b\rangle}_x})^{2\delta} = \dfrac{d\mu_x(a)d\mu_x(b)}{\big(g_{\infty, x}(a, b)\big)^{2\delta}}$$
where $a, b \in \partial^{\infty}(X)$, $x \in X$ and $g_{\infty, x}$ is the visual metric from Definition \ref{visualmetric}. 
One can check that due to the scale factor $(e^{{\langle a, b\rangle}_x})^{2\delta}$, the construction of $\nu_x$ is independent of choice of basepoint. We can then denote $\nu = \nu_x$ where $x$ is any arbitrary point in $X$. We can thus define a $\Gamma$-invariant Radon measure on $\partial^{\infty}(X) \times \partial^{\infty}(X) \setminus \Delta$, which we will also denote as $\nu$. 

While $\nu$ is defined on $\partial^{\infty}(X) \times \partial^{\infty}(X) \setminus \Delta$, there is a natural extension of $\nu$ to a measure $m$ on $SX = (\partial^{\infty}(X) \times \partial^{\infty}(X) \setminus \Delta) \times \mathbb{R}$: 
$$m = \nu \times dt$$
where $dt$ is the usual Lebesgue measure on $\mathbb{R}$. Note that $m$ descends to a measure, $m_{\Gamma}$, on $SX/\Gamma$.

Suppose that $\widetilde{X}$ is a (Gromov) hyperbolic covering space with $\Gamma$ its deck group. Unsurprisingly, for a Gromov hyperbolic metric space $X= \widetilde{X}/\Gamma$, there is a one-to-one correspondence between \textit{geodesic currents}, $\Gamma$-invariant Radon measures on $\partial^{\infty}(X) \times \partial^{\infty}(X) \setminus \Delta$, and Radon measures on $SX = S\widetilde{X}/\Gamma$ that are invariant under the geodesic flow (see Theorem 2.2 of \cite{kai}). Thus, since $\nu$ is $\Gamma$-invariant, $m_{\Gamma}$ is invariant under geodesic flow.

\begin{definition}[Bowen-Margulis measures] The geodesic flow-invariant measure $m_{\Gamma}$ defined above is a \textit{Bowen-Margulis measure} on $SX / \Gamma$. 
\end{definition}

\begin{remark} Note that $d\nu(a, b) = d\nu(b, a)$, so it follows that $\nu$ is invariant under the action of $\mathbb{Z}/2\mathbb{Z}$. The Bowen-Margulis measure thus provides a ``natural" example of a geodesic current on $X$. Its relationship with another naturally-arising current, the \textit{Liouville current}, is an interesting open question for spaces that are not symmetric. While we will not define Liouville currents in this paper, we refer the reader to Example 5.4 of \cite{con} for a definition suited to compact quotients of Fuchsian buildings which extends to the surface amalgam setting.
\end{remark}

\subsection{Proof of Proposition \ref{ergodicity}} 

We are now ready to state a key theorem from \cite{kai}: 

\begin{theorem}[\cite{kai}, Theorem 2.6] Let $\widetilde{X}$ be a (Gromov) hyperbolic covering space satisfying Assumptions \ref{U} and \ref{C}. Let $\mu_x$ be a Patterson-Sullivan measure on $\partial^{\infty}(X)$ used to construct the geodesic current $\nu$ from before and $m_{\Gamma}$ the corresponding Bowen-Margulis measure that is invariant under geodesic flow on $SX$. Then either:
\begin{enumerate}
    \item $\mu_x(\Lambda_{\Gamma}) = 1$ and the geodesic flow on $SX = S\widetilde{X}/\Gamma$ is ergodic with respect to $m_{\Gamma}$ or 
    \item $\mu_x(\Lambda_{\Gamma}) = 0$ and the geodesic flow on $SX$ is completely dissipative with respect to $m_{\Gamma}$. 
\end{enumerate}
\end{theorem}

In the case where $X$ is a P-manifold with a locally CAT(-1) metric, recall from Theorem \ref{qc} that $\text{supp}(\mu_x) = \Lambda_{\Gamma}$; therefore, the geodesic flow on SX is ergodic with respect to $m_{\Gamma}$, as desired. 

\subsubsection{An Application of Proposition \ref{ergodicity}} \label{birkhoffapp}

One application of the ergodicity of geodesic flow is there exists $v \in SX$ with dense orbit under the geodesic flow map. Indeed, recall the Birkhoff Ergodic Theorem:

\begin{theorem}[Birkhoff Ergodic Theorem]\label{birkhofftheorem} Let $(Y, \mathscr{B}, \mu)$ be a probability space, and let $T: Y \rightarrow Y$ be an ergodic measure-preserving transformation. Let $A \in \mathscr{B}$ be a measurable set of positive measure $\mu(A) > 0$. Then for all $f \in \mathscr{L}^1(Y, \mathscr{B}, \mu)$ and $\mu$-almost everywhere $y \in Y$: 
\begin{equation}
    \label{birkhoff}
    \lim\limits_{n \rightarrow \infty} \frac{1}{n}\sum\limits_{k = 0}^{n - 1} f(T^k(y)) = \int f d\mu.
\end{equation}

\end{theorem}

If $\Gamma$ is a nonelementary discrete group acting properly and cocompactly by isometries on a Gromov hyperbolic space $X$, then $m_{\Gamma}$ is finite (see Corollary 4.17 and Section 3.2 of \cite{coulon}) and in particular can be scaled to a probability measure, so we can apply the Birkhoff Ergodic Theorem. For example, if $X$ is a locally CAT(-1) P-manifold and $\Gamma = \pi_1(X)$, then we can apply the Birkhoff Ergodic Theorem. By the Hopf-Tsuji-Sullivan Theorem (Theorem 4.2) and Proposition 4.4 of \cite{coulon}, since the geodesic flow map is ergodic with respect to $m_{\Gamma}$, the flow map is conservative, so $m_{\Gamma}$ has full support on $SX/\Gamma$. 

In particular, if $f$ is the indicator function $\chi_{A}$, $\int f d\mu = \mu(A)$ while the left hand side of Equation \ref{birkhoff} denotes the frequency that the orbit of $T$ visits $A$. In our scenario, $T^k$ is the geodesic flow map $\phi_k$ on $SX$, which is invariant with respect to the scaled probability measure of $m_{\Gamma}$ with full support on $SX$. We thus conclude for every simple, thick locally CAT(-1) P-manifold, there exists a geodesic in $SX$ that intersects any open neighborhood of any $u \in SX$ (where $X$ satisfies the conditions of Theorem \ref{birkhofftheorem}), as every open set of $SX$ has nonzero measure. 

\printbibliography

\end{document}